\RequirePackage{fix-cm}
\documentclass{svjour3}                     
\smartqed  
\usepackage{graphicx}
\usepackage{amsmath} %
\usepackage{amsfonts} %
\usepackage{amssymb} %
\usepackage{mathtools}
\usepackage {verbatim}
\usepackage{csquotes}%
\usepackage{algorithm}
\usepackage{algorithmicx} 
\usepackage{algpseudocode}

\usepackage{float}
\usepackage{booktabs}
\usepackage{bookmark}
\usepackage{subcaption}
\usepackage{epstopdf}
\graphicspath{{figure_eps/}}
\begin{document}

\title{A Locally Divergence-Free Oscillation-Eliminating Discontinuous Galerkin Method for Ideal Magnetohydrodynamic Equations
}

\titlerunning{A LDF-OEDG Method for Ideal MHD Equations}        

\author{Wei Zeng         \and
        Qian Wang 
}


\institute{W. Zeng \at
              Mechanics Division, Beijing Computational Science Research Center, Beijing 100193, China \\
              \email{zengwei@csrc.ac.cn}           
           \and
           Q. Wang (Corresponding author) \at
              Mechanics Division, Beijing Computational Science Research Center, Beijing 100193, China\\
              \email{qian.wang@csrc.ac.cn}
}

\date{Received: date / Accepted: date}

\maketitle

\begin{abstract}
Numerical simulations of ideal compressible magnetohydrodynamic (MHD) equations are challenging, as the solutions are required to be magnetic divergence-free for general cases as well as oscillation-free for cases involving discontinuities.  To overcome these difficulties, we develop a locally divergence-free oscillation-eliminating discontinuous Galerkin (LDF-OEDG) method for ideal compressible MHD equations. In the LDF-OEDG method, the numerical solution is advanced in time by using a strong stability preserving Runge-Kutta scheme. 
Following the solution update in each Runge-Kutta stage, an oscillation-eliminating (OE) procedure is performed to suppress spurious oscillations near discontinuities by damping the modal coefficients of the numerical solution.  
Subsequently, on each element, the magnetic filed of the oscillation-free DG solution is projected onto a local divergence-free space, to satisfy the divergence-free condition. 
The OE procedure and the LDF projection are fully decoupled from the Runge-Kutta stage update, and can be non-intrusively integrated into existing DG codes as independent modules. 
The damping equation of the OE procedure can be solved exactly, making the LDF-OEDG method remain stable under normal CFL conditions. These features enable a straightforward implementation of a high-order LDF-OEDG solver, which can be used to efficiently simulate the ideal compressible MHD equations.
Numerical results for benchmark cases demonstrate the high-order accuracy, strong shock capturing capability and robustness of the LDF-OEDG method.

\end{abstract}

\section{Introduction}
\label{Sect1}
Magnetohydrodynamics (MHD) couples fluid dynamics with electromagnetism to describe the macroscopic behavior of conducting fluids such as plasma. Ideal MHD equations are applicable when relativistic, viscous and resistive effects can be neglected. The ideal MHD equations have been widely used in areas such as astrophysics, space physics and engineering. It is of great importance to design accurate and robust numerical methods for ideal MHD equations to efficiently simulate problems in these areas. There are two main difficulties in numerically solving the ideal compressible MHD equations, namely the magnetic divergence-free constraint for general cases, and the elimination of spurious oscillations for cases involving discontinuities.

The first difficulty comes from the physical principle of non-existence of magnetic monopoles. Based on this principle, the divergence-free (DF) constraint on the magnetic filed holds for the exact solution as long as it does initially. For a simulation, the numerical method must preserve the DF property, of which the violation may result in nonphysical solution or even break-down of the simulation. 
Tremendous efforts have been made in developing divergence-free numerical methods in the past few decades, among which the prominent ones are the constrained transport (CT) methods and the eight-wave (EW) methods. The former approach, named CT by Evans and Hawley \cite{Evans1988Simulation}, refers to a specific finite difference discretization on a staggered grid that ensures the maintenance of magnetic divergence. Various adaptations have been developed within different frameworks, see \cite{DeVore1991Flux,Dai1998A,Londrillo2004On,Gardiner2005An,Xu2016Divergence}. Additionally, unstaggered CT methods have also been devised (see e.g., \cite{Rossmanith2006An,Helzel2011An,Mishra2012Constraint,Helzel2013A,Christlieb2014Finite}), typically based on numerically evolving the magnetic potential. 
The latter approach refers to the EW formulation \cite{Gombosi1994Axisymmetric} that is based on proper discretization of the Godunov form \cite{Godunov1972Symmetric} of ideal MHD equations. This approach incorporates an additional source term proportional to magnetic divergence and a simple modification of Riemann solver. 
 It is noted that the Godunov form of ideal MHD equations has its advantage in terms of positivity, see \cite{Wu2018A,Wu2019Provably,Wu2021Provably,Liu2021A,Zhang2021A} for specific applications. 
 Besides the CT methods and the EW methods, existing divergence-free methods can be roughly divided into four categories: hyperbolic divergence cleaning \cite{Dedner2002Hyperbolic}, projection \cite{Brackbill1980The}, globally divergence-free (GDF) \cite{Li2011Central,Li2012Arbitrary,Fu2018Globally}, and locally divergence-free (LDF) \cite{Cockburn2004Locally,Li2005Locally,Yakovlev2013Locally}. 
 In the LDF approach, the magnetic field is approximated in a local divergence-free polynomial space, so that the divergence-free condition can be automatically satisfied inside each element. 
 The LDF solution can be also obtained by projecting the solution computed in the standard polynomial space onto a local divergence-free polynomial space \cite{Liu2019Locally}. The projection does not interfere with the spatial discretization, thus enabling a non-intrusive implementation of the LDF approach.

The second difficulty lies in eliminating spurious oscillations near discontinuities while preserving accuracy in smooth regions. To the best of our knowledge, four kinds of discontinuity capturing approaches have been investigated in previous studies. Among these, slope limiter is a prominent technique, which is used to regularize numerical solutions. With proper adjustments of parameters, the limiters proved to work well, as noted in \cite{Cockburn1998The,Zhong2013A}. In particular, some classical and central discontinuous Galerkin (DG) methods are paired with the minmod total variation bounded (TVB) slope limiter to solve MHD equations, see, for example, \cite{Li2011Central,Li2012Arbitrary,Fu2018Globally,Cockburn2004Locally,Li2005Locally,Yakovlev2013Locally}. Another well-regarded strategy involves incorporating artificial diffusion terms into the weak formulations, ensuring that the adjusted formulation aligns with the original system as the mesh size diminishes \cite{Barter2010Shock,Hiltebrand2014Entropy,Xin2006Viscous}. Moreover, spectral filters constitute a third regularization method, prized for its simple implementation and low computational complexity. The spectral viscosity approach is a notable example \cite{Maday1993Legendre}. In \cite{Hesthaven2007nodal} and \cite{Meister2012Application}, the exponential and adaptive filters are applied to DG methods, respectively. Additionally, the damping technique recently developed in \cite{Lu2021An,Liu2022An,Peng2023OEDG}, can be regarded as the fourth approach, which is inspired by the spectral filters. In the DG setting, the oscillation-free discontinuous Galerkin (OFDG) method \cite{Lu2021An,Liu2022An} incorporates a damping term to control oscillations, rather than a post-processor like slope limiter. Unfortunately, for discontinuous problems, the damping term renders the semi-discrete OFDG system highly stiff, imposing severe restrictions on step sizes for explicit time integration to be stable. Recently, on the basis of OFDG method, Peng et al. \cite{Peng2023OEDG} developed a simple, efficient and non-intrusive oscillation eliminating (OE) procedure for the DG method. In the OE procedure, an ordinary differential equation (ODE) is solved exactly to damp the modal coefficients of the numerical solution. The numerical results for compressible Euler equations in \cite{Peng2023OEDG} demonstrate that the oscillation-eliminating discontinuous Galerkin (OEDG) method is able to suppress spurious oscillations near discontinuities while preserving accuracy in smooth regions, and also stable under normal CFL conditions. Besides, the OE procedure is free of characteristic decomposition. Based on these features, the OEDG method has great potential in more complex applications such as the ideal MHD equations.

In this work, we develop a locally divergence-free oscillation-eliminating discontinuous Galerkin (LDF-OEDG) method for ideal compressible MHD equations to overcome the aforementioned two difficulties. In the LDF-OEDG method, the LDF projection and OE procedure are used to obtain divergence- and oscillation-free solutions. 
Specifically, all the conservative variables are approximated in the standard piece-wise polynomial space. The semi-discrete DG scheme is integrated in time by using a strong stability preserving (SSP) Runge-Kutta scheme \cite{Shu1988efficient} to update the degrees of freedom. Following the solution update in each Runge-Kutta stage, an OE procedure is performed to suppress spurious oscillations near discontinuities, by applying an exact damping operator to the modal coefficients of the DG solution. Subsequently, a LDF projection is performed to enforce the divergence-free constraint on the numerical solution. On each element, the OEDG solution is projected onto a local divergence-free polynomial space. Both the OE procedure and the LDF projection are non-intrusive, in the sense that they are fully decoupled from the Runge-Kutta stage update. Therefore, they are two post-processing procedures that can be seamlessly integrated into existing Runge-Kutta discontinuous Galerkin (RKDG) codes as independent modules, thus enabling a straightforward and efficient implementation of a LDF-OEDG solver. 
The LDF-OEDG method is applied to a series of two-dimensional benchmark ideal compressible MHD test cases. Numerical results demonstrate the high-order accuracy, strong shock capturing capability and robustness of the LDF-OEDG method.

The remainder of this paper is organized as follows. Section \ref{Sect2} presents the locally divergence-free Runge-Kutta discontinuous Galerkin method for ideal compressible MHD equations. Section \ref{Sect3} presents the oscillation-eliminating procedure for the locally divergence-free discontinuous Galerkin method. Numerical results are presented in Section \ref{Sect4} and concluding remarks are given in Section \ref{Sect5}.

\section{Locally divergence-free discontinuous Galerkin method for ideal compressible MHD equations}
\label{Sect2}
\subsection{Governing equations} \label{subsec:governing_equations}
In this paper, we consider the ideal compressible MHD equations that can be formulated in a conservative form as
\begin{equation}\label{MHD}	 
	\partial_{t} \mathbf{U} + \nabla \cdot \mathbf{F}\left(\mathbf{U}\right)= 0, 
\end{equation}
where $\mathbf{U}$ is the vector of conservative variables and $\mathbf{F}$ is the vector of convective fluxes defined by
\begin{equation*}\label{MHD_U_F}
	\mathbf{U}= \begin{pmatrix}
		\rho\\ 	\rho \mathbf{u}\\ 	E\\ 	\mathbf{B} 
	\end{pmatrix}, \quad
	\mathbf{F}= \begin{pmatrix}
		\rho \mathbf{u}\\ 	\rho \mathbf{u}\otimes \mathbf{u}+p_{tot}\mathcal{I}-\mathbf{B}\otimes \mathbf{B}\\ 	\mathbf{u}\left( E +p_{tot} \right) -\mathbf{B}\left( \mathbf{u}\cdot \mathbf{B} \right)\\ 	\mathbf{u}\otimes \mathbf{B}-\mathbf{B}\otimes \mathbf{u}
	\end{pmatrix}.
\end{equation*}
Here $\rho$ is the density, $\mathbf{u} \in \mathbb{R}^{d}$ is the velocity, $\mathbf{B} \in \mathbb{R}^{d}$ is the magnetic field, $p_{tot}=p+\frac{\|\mathbf B\|^2}{2}$ is the total pressure, $p$ is the thermal pressure, $E$ is the total energy, $\mathcal{I}$ is the $d\times d$ identity matrix, and $d$ is the spatial dimensionality. The total energy is computed by
\begin{equation*}
	E=\dfrac{p}{\gamma-1}+\frac{1}{2}\rho \left\|\mathbf u\right\|^2+\frac{1}{2}\left\|\mathbf B\right\|^2,
\end{equation*}
where $\gamma$ is the adiabatic index.

Based on the physical principle of non-existence of magnetic monopoles, the magnetic field $\mathbf B$ should satisfy a divergence-free (DF) constraint
\begin{equation}\label{divB0} 
	\nabla\cdot\mathbf B=0.
\end{equation}
In fact, according to the governing equation of $\mathbf{B}$ in \eqref{MHD}
\begin{equation*}
	\frac{\partial \mathbf B}{\partial t} + \nabla\cdot(\mathbf u\otimes \mathbf B - \mathbf B\otimes\mathbf u) = 0, 
\end{equation*}
we have
\begin{equation*}
	\displaystyle\frac{\partial(\nabla\cdot\mathbf B)}{\partial t}=\nabla\cdot\frac{\partial \mathbf B}{\partial t}=\nabla\cdot(\nabla\cdot(\mathbf B\otimes\mathbf u-\mathbf u\otimes\mathbf B)) = \nabla\cdot(\nabla \times \left( \mathbf{u}\times \mathbf{B} \right))=0,
\end{equation*}
that is,
\begin{equation*}
	\nabla\cdot\mathbf B(\boldsymbol{x},t) = \nabla\cdot\mathbf B(\boldsymbol{x},0),\ \forall~t>0, 
\end{equation*}
which implies that the DF constraint \eqref{divB0} holds for the exact solution if the initial divergence of $\mathbf B$ is zero.

\subsection{Runge-Kutta discontinuous Galerkin method}\label{subsec:rkdg}

For the sake of presentation, we consider the two-dimensional case ($d=2$) of \eqref{MHD} to illustrate the spatial and temporal discretizations. A two-dimensional computational domain $\Omega$ is partitioned into non-overlapping triangular/rectangular control volumes. Let $\mathcal{X}_h$ be a partition of $\Omega$. 
On each control volume, the distribution of each conservative variable is approximated by a polynomial. The piece-wise polynomial distribution of the numerical solution $\mathbf{U}_h$ is discontinuous across cell interfaces. Specifically, on a control volume $K \in \mathcal{X}_h$, solution $\mathbf{U}$ is approximated by a polynomial of degree $k$ as 
\begin{equation}\label{Uh}
	\mathbf{U}_h \left(\boldsymbol{x},t\right)= \sum_{\left| \boldsymbol{\alpha} \right|=0}^{k} \mathbf{U}^{\left(\boldsymbol{\alpha}\right)}_K \left(t\right) \phi^{\left(\boldsymbol{\alpha}\right)}_K \left(\boldsymbol{x}\right), \quad \forall~\boldsymbol{x} \in K,
\end{equation}
where $\left\{\phi^{\left(\boldsymbol{\alpha}\right)}_K\right\}^k_{\left| \boldsymbol{\alpha} \right|=0}$ is a polynomial basis and $\boldsymbol{\alpha}=\left(\alpha_1,\alpha_2\right)$ is the multi-index vector with $\left| \boldsymbol{\alpha} \right|=\alpha_1 + \alpha_2$. The total number of basis functions is $\left(k+1\right)\left(k+2\right)/2$. The basis coefficients $\left\{\mathbf{U}^{\left(\boldsymbol{\alpha}\right)}_K \right\}^k_{\left| \boldsymbol{\alpha} \right|=0}$ are the degrees of freedom (DOF).

The semi-discrete discontinuous Galerkin (DG) scheme for the governing equation \eqref{MHD} reads
\begin{equation}\label{semi-DG}
	\int_K \left(\mathbf{U}_h \right)_t \phi^{\left(\boldsymbol{\alpha}\right)}_K \ d\boldsymbol{x} + \oint_{\partial K} \mathbf{\hat F} \cdot \mathbf{n} \ \phi^{\left(\boldsymbol{\alpha}\right)}_K \ ds - \int_K \mathbf{F} \cdot \nabla \phi^{\left(\boldsymbol{\alpha}\right)}_K \ d\boldsymbol{x} =0, \quad \left| \boldsymbol{\alpha} \right| \le k,
\end{equation}
where $\partial K$ is the boundary of $K$, $\mathbf{n}$ is the outward unit normal of $\partial K$, and $\mathbf{\hat F}$ is a numerical flux on $\partial K$. Gauss quadrature rules are used to compute both the surface and volume integrals in \eqref{semi-DG}. As $\mathbf{U}_h$ is discontinuous across cell interfaces, $\mathbf{\hat F} \cdot \mathbf{n}$ is computed by using an approximate Riemann solver. In this work, the local Lax-Friedrichs flux
\begin{equation*}
	\mathbf{\hat F} \cdot \mathbf{n} = \dfrac{1}{2} \left[\mathbf{F}\left(\mathbf{U}_L\right) \cdot \mathbf{n} + \mathbf{F}\left(\mathbf{U}_R\right) \cdot \mathbf{n} \right] - \dfrac{1}{2} \lambda_{max} \left(\mathbf{U}_R - \mathbf{U}_L\right),
\end{equation*}
is computed on each surface quadrature point $\boldsymbol{x}_g \in \partial K$, where $\mathbf{U}_L$ is the left state (inside $K$), $\mathbf{U}_R$ is the right state (outside $K$), and $\lambda_{max}$ is an estimate of the local maximum wave speed in direction $\mathbf{n}$.
In this work, an orthogonal basis $\left\{\phi^{\left(\boldsymbol{\alpha}\right)}_K\right\}^k_{\left| \boldsymbol{\alpha} \right|=0}$ which has the property 
\begin{equation*}
	\int_K   \phi^{\left(\boldsymbol{\alpha}_i\right)}_K  \phi^{\left(\boldsymbol{\alpha}_j\right)}_K \ d\boldsymbol{x} = \delta_{ij} \int_K   \left(\phi^{\left(\boldsymbol{\alpha}_i\right)}_K \right)^2 \ d\boldsymbol{x},
\end{equation*}
is used to reduce the semi-discrete DG scheme \eqref{semi-DG} to 
\begin{equation*}\label{orthogonal-semi-DG}
	\dfrac{d\mathbf{U}_K^{\left(\boldsymbol{\alpha}\right)}}{dt}  \int_K  \left(\phi^{\left(\boldsymbol{\alpha}\right)}_K \right)^2 \ d\boldsymbol{x} + \oint_{\partial K} \mathbf{\hat F} \cdot \mathbf{n} \ \phi^{\left(\boldsymbol{\alpha}\right)}_K \ ds - \int_K \mathbf{F} \cdot \nabla \phi^{\left(\boldsymbol{\alpha}\right)}_K \ d\boldsymbol{x} =0, \quad \left|\boldsymbol{\alpha} \right| \le k.
\end{equation*}
The orthogonal basis can be obtained by applying the Gram-Schmidt process to a Taylor basis \cite{Hesthaven2007nodal,Zhu2008runge,Li2014multi}. For instance, on a rectangular element $K=\left[x_{1}, x_{2}\right] \times \left[y_{1}, y_{2}\right]$, an orthogonal polynomial basis of degree $k=2$ is
\begin{equation*}\label{orthogonal-basis}
	\begin{aligned}
		\phi^{\left(\left(0,0\right)\right)}_K &= 1, \\
		\phi^{\left(\left(1,0\right)\right)}_K &= X,\ \phi^{\left(\left(0,1\right)\right)}_K= Y, \\
		\phi^{\left(\left(2,0\right)\right)}_K &= X^2-\frac{1}{3},\ \phi^{\left(\left(1,1\right)\right)}_K= XY, \ \phi^{\left(\left(0,2\right)\right)}_K= Y^2-\frac{1}{3},
	\end{aligned}
\end{equation*}
where 
\begin{equation*}
	X=\frac{x-x_K}{h_x/2}, Y=\frac{y-y_K}{h_y/2},  h_x= x_{2}-x_{1}, h_y=y_{2}-y_{1}, x_K= \frac{x_{1}+x_{2}}{2}, y_K=\frac{y_{1}+y_{2}}{2}.
\end{equation*}

The semi-discrete DG scheme \eqref{semi-DG} can be rewritten in an ODE form
\begin{equation}\label{semi-DG-ODE}
	\dfrac{d}{dt} \mathbf{U}_h = \mathcal T_f \left(\mathbf{U}_h\right),
\end{equation} 
which can be integrated in time to update the DG solution $\mathbf{U}_h$ in a step-by-step manner. 
In this work, a third-order strong-stability-preserving Runge-Kutta (SSPRK3) method \cite{Shu1988efficient} is used for time integration. Specifically, for time step $n$, the solution at the next time step is computed by
\begin{equation}\label{SSP-RK}
	\begin{aligned}
		\mathbf U_h^{n,1} & =\mathbf U_h^n+\tau \mathcal T_f\left(\mathbf U_h^n\right), \\
		\mathbf U_h^{n,2} & =\frac{3}{4} \mathbf U_h^n+\frac{1}{4}\left(\mathbf U_h^{n,1}+\tau \mathcal T_f\left(\mathbf U_h^{n,1}\right)\right), \\
		\mathbf U_h^{n+1} & =\frac{1}{3} \mathbf U_h^n+\frac{2}{3}\left(\mathbf U_h^{n,2}+\tau \mathcal T_f\left(\mathbf U_h^{n,2}\right)\right),
	\end{aligned}
\end{equation}
with $\tau= t_{n+1}  - t_n$ being the time step size.

\subsection{Locally divergence-free Runge-Kutta discontinuous Galerkin method}\label{subsec:ldf-dg}

The RKDG method described in Section \ref{subsec:rkdg} can be directly applied to the mass, momentum and energy conservation equations in \eqref{MHD}. However, modifications are needed for the RKDG method to satisfy the divergence-free constraint \eqref{divB0} when solving the governing equations of the magnetic filed $\mathbf{B}$.

A direct approach to obtain magnetic divergence-free solutions is to use a locally divergence-free (LDF) basis that has the property
\begin{equation*}
	\nabla \cdot \boldsymbol{\psi}_K \left(\boldsymbol{x}\right)=0, \quad \forall \ \boldsymbol{x} \in K.
\end{equation*}
Therefore, as a linear combination of the LDF basis functions, the magnetic filed of the LDF-DG solution
\begin{equation*}\label{LDF-B}
	\mathbf{B}_h \left(\boldsymbol{x},t\right)= \sum_{l} B^{\left(l\right)}_K \left(t\right)\boldsymbol{\psi}_K^{\left(l\right)} \left(\boldsymbol{x}\right), \quad \forall \ \boldsymbol{x} \in K,
\end{equation*}
can automatically satisfy the divergence-free constraint. For a rectangular element $K=\left[x_{1}, x_{2}\right] \times \left[y_{1}, y_{2}\right]$, a valid orthogonal divergence-free basis of degree $k=2$ is
\begin{equation*}
	\begin{aligned}
		&\boldsymbol{\psi}^{\left(1\right)}_K = \begin{pmatrix}
			\phi^{\left(\left(0,0\right)\right)}_K \\
			0
		\end{pmatrix},
		\boldsymbol{\psi}^{\left(2\right)}_K = \begin{pmatrix}
			0\\
			\phi^{\left(\left(0,0\right)\right)}_K
		\end{pmatrix},
		\\
		&\boldsymbol{\psi}^{\left(3\right)}_K = \begin{pmatrix}
			h_x \phi^{\left(\left(1,0\right)\right)}_K \\
			-h_y \phi^{\left(\left(0,1\right)\right)}_K
		\end{pmatrix},
		\boldsymbol{\psi}^{\left(4\right)}_K = \begin{pmatrix}
			\phi^{\left(\left(0,1\right)\right)}_K\\
			0
		\end{pmatrix}, 
		\boldsymbol{\psi}^{\left(5\right)}_K = \begin{pmatrix}
			0 \\
			\phi^{\left(\left(1,0\right)\right)}_K
		\end{pmatrix},
		\\
		&\boldsymbol{\psi}^{\left(6\right)}_K = \begin{pmatrix}
			h_x \phi^{\left(\left(2,0\right)\right)}_K \\
			-2h_y \phi^{\left(\left(1,1\right)\right)}_K 
		\end{pmatrix},
		\boldsymbol{\psi}^{\left(7\right)}_K = \begin{pmatrix}
			2h_x \phi^{\left(\left(1,1\right)\right)}_K \\
			-h_y \phi^{\left(\left(0,2\right)\right)}_K 
		\end{pmatrix}, 
		\boldsymbol{\psi}^{\left(8\right)}_K = \begin{pmatrix}
			\phi^{\left(\left(0,2\right)\right)}_K \\
			0
		\end{pmatrix}, 
		\boldsymbol{\psi}^{\left(9\right)}_K = \begin{pmatrix}
			0 \\
			\phi^{\left(\left(2,0\right)\right)}_K 
		\end{pmatrix},
	\end{aligned}
\end{equation*}
which spans a space $\Psi_K= \left\{\boldsymbol{\psi}^{\left(1\right)}_K, \cdots, \boldsymbol{\psi}^{\left(9\right)}_K\right\} \in \mathbb{R}^{2 \times 9}$. $\Psi_K$ is a subspace of the space spanned by the full polynomial basis, i.e.,
\begin{equation*}
	\Psi_K \subset \Phi_K=\left\{\begin{pmatrix}
		\phi^{\left( \left(0,0\right) \right)}_K \\
		0
	\end{pmatrix}, \cdots, \begin{pmatrix}
		\phi^{\left( \left(0,2\right) \right)}_K \\
		0
	\end{pmatrix}, 
	\begin{pmatrix}
		0 \\
		\phi^{\left( \left(0,0\right) \right)}_K
	\end{pmatrix}, \cdots, \begin{pmatrix}
		0 \\
		\phi^{\left( \left(0,2\right) \right)}_K
	\end{pmatrix}\
	\right\} \in \mathbb{R}^{2 \times 12}.
\end{equation*}
The semi-discrete LDF-DG scheme of the magnetic conservation equation is
\begin{equation}\label{LDF-semi-DG}
	\int_K \left(\mathbf{B}_h \right)_t \cdot \boldsymbol{\psi} \ d\boldsymbol{x} + \oint_{\partial K} \left( \mathbf{\hat F}_{\mathbf{B}} \cdot \mathbf{n} \right) \cdot \boldsymbol{\psi} \ ds - \int_K \mathbf{F}_{\mathbf{B}} \cdot \nabla \boldsymbol{\psi} \ d\boldsymbol{x} =0, \quad \forall ~ \boldsymbol{\psi} \in \Psi_K,
\end{equation}
which is coupled with the mass, momentum and energy parts of usual DG scheme \eqref{semi-DG} to solve the complete ideal compressible MHD equations. A LDF-RKDG scheme can be obtained by applying the SSP Runge-Kutta time integration \eqref{SSP-RK} to \eqref{LDF-semi-DG}.

The LDF-RKDG approach is equivalent to the classical RKDG approach followed by a projection, which is derived as follows. As the SSP Runge-Kutta scheme is a convex combination of forward Euler schemes, we only need to prove the equivalence based on forward Euler time integration. The divergence-free magnetic filed at time step $n$ is denoted as $\mathbf{B}^n_h$. In LDF-DG approach, the solution is updated by
\begin{equation}\label{Euler-LDF-DG}
	\begin{aligned}
	\int_K \mathbf{B}^{n+1}_h  \cdot \boldsymbol{\psi} \ d\boldsymbol{x} =& \int_K \mathbf{B}^{n}_h \cdot \boldsymbol{\psi} \ d\boldsymbol{x}\\
	  &- \tau \left(  \oint_{\partial K} \left( \mathbf{\hat F}_{\mathbf{B}} \cdot \mathbf{n} \right) \cdot \boldsymbol{\psi} \ ds - \int_K \mathbf{F}_{\mathbf{B}} \cdot \nabla \boldsymbol{\psi} \ d\boldsymbol{x} \right), \ \forall ~ \boldsymbol{\psi} \in \Psi_K.
	\end{aligned}
\end{equation}
In classical DG approach, the magnetic filed $\mathbf{B}^{n+1,*}_h$ is computed by
\begin{equation*}\label{Euler-DG}
	\begin{aligned}
	\int_K \mathbf{B}^{n+1,*}_h  \cdot \boldsymbol{\psi} \ d\boldsymbol{x} =& \int_K \mathbf{B}^{n}_h \cdot \boldsymbol{\psi} \ d\boldsymbol{x}\\
	  &- \tau \left(  \oint_{\partial K} \left( \mathbf{\hat F}_{\mathbf{B}} \cdot \mathbf{n} \right) \cdot \boldsymbol{\psi} \ ds - \int_K \mathbf{F}_{\mathbf{B}} \cdot \nabla \boldsymbol{\psi} \ d\boldsymbol{x} \right), \ \forall ~ \boldsymbol{\psi} \in \Phi_K.
	\end{aligned}
\end{equation*}
As $\Psi_K \subset \Phi_K$, we have
\begin{equation}\label{DG-projection-LDF-space}
	\begin{aligned}
	\int_K \mathbf{B}^{n+1,*}_h  \cdot \boldsymbol{\psi} \ d\boldsymbol{x} =& \int_K \mathbf{B}^{n}_h \cdot \boldsymbol{\psi} \ d\boldsymbol{x}\\
	  &- \tau \left(  \oint_{\partial K} \left( \mathbf{\hat F}_{\mathbf{B}} \cdot \mathbf{n} \right) \cdot \boldsymbol{\psi} \ ds - \int_K \mathbf{F}_{\mathbf{B}} \cdot \nabla \boldsymbol{\psi} \ d\boldsymbol{x} \right), \ \forall ~ \boldsymbol{\psi} \in \Psi_K.
	\end{aligned}
\end{equation}
By comparing \eqref{Euler-LDF-DG} and \eqref{DG-projection-LDF-space}, we obtain the following relation
\begin{equation}\label{Projection}
	\int_K \mathbf{B}^{n+1}_h \cdot \boldsymbol{\psi} \ d\boldsymbol{x} = \int_K \mathbf{B}^{n+1,*}_h  \cdot \boldsymbol{\psi} \ d\boldsymbol{x}, \quad \forall ~ \boldsymbol{\psi} \in \Psi_K,
\end{equation}
which implies that the LDF-DG solution $\mathbf{B}^{n+1}_h$ is the projection of the classical DG solution $\mathbf{B}^{n+1,*}_h$ onto the divergence-free subspace $\Psi_K$. Therefore, the LDF-DG scheme is equivalent to the classical DG scheme followed by a projection. For the rectangular element $K$ using polynomial basis of degree $k=2$, the divergence-free solution is computed by
\begin{equation}\label{formula-projection}
	\mathbf{B}^{n+1}_h= \sum_{l=1}^{9} B^{\left(l\right)}_K \boldsymbol{\psi}^{\left(l\right)}_K, \quad B^{\left(l\right)}_K= \dfrac{\int_K \mathbf{B}^{n+1,*}_h \cdot \boldsymbol{\psi}^{\left(l\right)}_K \ d\boldsymbol{x}}{\int_K \boldsymbol{\psi}^{\left(l\right)}_K \cdot \boldsymbol{\psi}^{\left(l\right)}_K \ d\boldsymbol{x}}.
\end{equation}

The equivalence relation \eqref{Projection} can significantly simplify the implementation of LDF-RKDG method, in the sense that one only needs to perform a classical RKDG method with a projection step at the end of each Runge-Kutta stage. State-of-the-art techniques for the RKDG method such as the oscillation-eliminating procedure that will be presented in Section \ref{Sect3}, can be directly extended to the LDF-RKDG method, as the LDF projection is non-intrusive and does not interfere with the DG spatial discretization.

\section{Oscillation-eliminating procedure}\label{Sect3}

The ideal compressible MHD equations \eqref{MHD} are nonlinear hyperbolic conservation laws, which may develop discontinuities within finite time even starting from smooth initial conditions. The high-order accurate LDF-RKDG scheme presented in Section \ref{Sect2} produces spurious oscillations near discontinuities, which may lead to nonphysical solutions or numerical instabilities. Suppressing oscillations in the vicinity of discontinuities while preserving high-order accuracy in smooth regions is a significant challenge in numerically solving the ideal compressible MHD equations. To overcome this difficulty, we develop a LDF-OEDG method that incorporates the recently developed oscillation-eliminating (OE) procedure \cite{Peng2023OEDG} into the LDF-RKDG framework. The OE procedure controls oscillations by damping modal coefficients of the numerical solution. The damping ODE of the OE procedure can be solved exactly, making the LDF-OEDG scheme remain stable under normal CFL conditions. The OE procedure does not affect the order of accuracy, as the damping term is high-order accurate in smooth regions. Furthermore, the OE procedure is non-intrusive in the sense that it is fully decoupled from the Runge-Kutta stage update as well as the LDF projection, enabling a straightforward implementation of the LDF-OEDG method. Details of the OE procedure will be presented in the remaining of this section.

\subsection{Oscillation-eliminating discontinuous Galerkin method} \label{subsec:OEDG}
The semi-discrete OEDG scheme for ideal compressible MHD equations \eqref{MHD} is obtained through an operator-splitting of the following damped semi-discrete DG scheme
\begin{equation}\label{damped-dg}
	\begin{aligned}
	\int_K \left(\mathbf{U}_h \right)_t \phi^{\left(\boldsymbol{\alpha}\right)}_K \ d\boldsymbol{x} =& -\oint_{\partial K} \mathbf{\hat F} \cdot \mathbf{n} \ \phi^{\left(\boldsymbol{\alpha}\right)}_K \ ds \ + \int_K \mathbf{F} \cdot \nabla \phi^{\left(\boldsymbol{\alpha}\right)}_K \ d\boldsymbol{x} \\
	 &- \sum_{m=0}^k \delta_K^m\left(\mathbf{U}_h\right) \int_K\left(\mathbf{U}_h-P^{m-1} \mathbf{U}_h\right) \phi^{\left(\boldsymbol{\alpha}\right)}_K \ d \boldsymbol{x},
	\end{aligned}
\end{equation}
in which the last term is a damping term that is added to suppress spurious oscillations near discontinuities. $P^{m}$ is the standard $L^2$ projection operator onto $\mathbb{P}^{m}\left(K\right)$, the space of polynomials of degree less than or equal to $m$ on $K$, for $m\ge0$. For any function $w$,  $P^{m} w \in \mathbb{P}^m \left(K\right)$ and 
\begin{equation*}
	\int_K  \left( P^{m} w -w\right) v \ d\boldsymbol{x}=0, \quad \forall \ v \in \mathbb{P}^m \left(K\right).
\end{equation*}
It is defined that $P^{-1}= P^0$. The damping term penalizes the deviation between the DG solution and its lower-order polynomial projections. The penalization effect is controlled by damping coefficients $\left\{ \delta_K^m\right\}^k_{m=0}$, which should be small in smooth regions and large near discontinuities. The damping coefficients should also be scale-invariant and evolution-invariant. The detailed formulae for the damping coefficients will be given in Section \ref{subsec:damping_coefficients}.

For discontinuous problems, the damped ODE \eqref{damped-dg} is highly stiff and unstable under normal CFL conditions. To reduce the stiffness of the system, Peng et. al. \cite{Peng2023OEDG} propose to split \eqref{damped-dg} into two ODEs, namely a conventional semi-discrete DG scheme \eqref{semi-DG-ODE} and a damping ODE, as
\begin{subequations}\label{semi-OEDG}
	\begin{align}
		&\dfrac{d}{dt} \mathbf{U}_h  = \mathcal T_f \left(\mathbf{U}_h\right), \label{normal-dg}\\
		& \dfrac{d}{dt} \mathbf{U}_{\sigma} = -\Sigma \left(\mathbf{U}_h\right) \mathbf{U}_{\sigma}, \label{damping-ode}
	\end{align}
\end{subequations}
where $\Sigma \left(\mathbf{U}_h\right)$ is a damping operator. Equation \eqref{damping-ode} corresponds to the following initial value problem
\begin{equation} \label{damping-problem}
	\begin{dcases}
		&\dfrac{d}{d\hat{t}} \int_K \mathbf{U}_{\sigma} \phi^{\left(\boldsymbol{\alpha}\right)}_K \ d\boldsymbol{x} = - \sum_{m=0}^k \delta_K^m\left(\mathbf{U}_h\right) \int_K\left(\mathbf{U}_\sigma-P^{m-1} \mathbf{U}_\sigma\right) \phi^{\left(\boldsymbol{\alpha}\right)}_K \ d \boldsymbol{x},\\
		&\mathbf{U}_{\sigma} |_{\hat{t}=0} = \mathbf{U}_h,
	\end{dcases}
\end{equation}
where $\hat{t}$ is a pseudo-time different from $t$. We define $\mathcal{F}_\tau$ as the solution operator of problem \eqref{damping-problem} with a pseudo-time step size $\tau$, i.e., $\left(\mathcal{F}_\tau \mathbf U_h \right) (\boldsymbol{x})= \mathbf{U}_\sigma(\boldsymbol{x},\tau)$. 
As pointed out in \cite{Peng2023OEDG}, with orthogonal DG basis functions $\left\{\phi^{\left(\boldsymbol{\alpha}\right)}_K \right\}^k_{\left|\boldsymbol{\alpha}\right|=0}$, the damping operator $\Sigma \left(\mathbf{U}_h\right)$ simplifies to a diagonal matrix, thus resulting in an exact solution
\begin{equation}\label{OE-solution}
	\mathcal{F}_\tau \mathbf{U}_h=\mathbf{U}_K^{(\mathbf{0})} \phi_K^{(\mathbf{0})}(\boldsymbol{x})+\sum_{j=1}^k \mathrm{e}^{-\tau \sum_{m=0}^j \delta_K^m\left(\mathbf{U}_h\right)} \sum_{|\boldsymbol{\alpha}|=j} \mathbf{U}_K^{(\boldsymbol{\alpha})} \phi_K^{(\boldsymbol{\alpha})}(\boldsymbol{x}),
\end{equation}
where the basis coefficients are $\mathbf{U}_K^{(\boldsymbol{\alpha})}=\int_K \mathbf{U}_h \phi_K^{(\boldsymbol{\alpha})} \ d\boldsymbol{x} /\|\phi_K^{(\boldsymbol{\alpha})}\|_{L^2(K)}^2$. It is shown in \eqref{OE-solution} that, the OE procedure does not affect the cell-average $\mathbf{\overline{U}}_K= \mathbf{U}_K^{(\mathbf{0})} \phi_K^{(\mathbf{0})}$, thus preserving the conservation property of the original DG scheme. The OE procedure actually damps the coefficients of high-order terms in numerical solution \eqref{Uh}.
The OEDG method based on the third-order SSP Runge-Kutta time stepping can be expressed as
\begin{equation*}
	\begin{aligned}
		\mathbf U_h^{n,1} & =\mathbf U_\sigma^n+\tau \mathcal T_f\left(\mathbf U_\sigma^n\right), & \mathbf U_\sigma^{n,1} & =\mathcal{F}_\tau \mathbf U_h^{n,1}, \\
		\mathbf U_h^{n,2} & =\frac{3}{4} \mathbf U_\sigma^n+\frac{1}{4}\left(\mathbf U_\sigma^{n,1}+\tau \mathcal T_f\left(\mathbf U_\sigma^{n,1}\right)\right), & \mathbf U_\sigma^{n,2} & =\mathcal{F}_\tau \mathbf U_h^{n,2}, \\
		\mathbf U_h^{n,3} & =\frac{1}{3} \mathbf U_\sigma^n+\frac{2}{3}\left(\mathbf U_\sigma^{n,2}+\tau \mathcal T_f\left(\mathbf U_\sigma^{n,2}\right)\right), & \mathbf U_\sigma^{n+1} & =\mathcal{F}_\tau \mathbf U_h^{n,3} .
	\end{aligned}
\end{equation*}
Thanks to the exact OE solver \eqref{OE-solution}, the OEDG method remains stable under normal CFL conditions. Furthermore, the OE procedure is fully decoupled from the Runge-Kutta stage update, thus can be implemented as an independent module that can be easily incorporated into existing DG codes. 

\subsection{Damping coefficients} \label{subsec:damping_coefficients}

There are several requirements for the design of the OE damping operator. First, the damping coefficients should be sufficiently large in the vicinity of discontinuities to suppress oscillations, while sufficiently small in smooth regions to preserve accuracy. Second, the damping coefficients should be scale- and evolution-invariant, to make the OE procedure perform consistently well for problems across various scales and wave speeds. We adopt the formulae in \cite{Peng2023OEDG} that generate damping coefficients satisfying the aforementioned properties.

The damping coefficient $\delta_K^m\left(\mathbf{U}_h\right)$ is computed as
\begin{equation}\label{delta}
	\delta_K^m \left(\mathbf{U}_h\right)= \sum_{e \in \partial K} \beta_e \dfrac{\sigma^m_{e, K} \left(\mathbf{U}_h\right)}{h_{e, K}},
\end{equation}
where $h_{e, K}= \sup_{\boldsymbol{x} \in K} \mathrm{dist} \left(\boldsymbol{x},e\right)$ is a characteristic length scale for cell interface $e$, and $\beta_e$ is the local maximum wave speed in the normal direction $\mathbf{n}_e$ evaluated at the cell-average $\mathbf{\overline{U}}_K$. Let $u^{\left(i\right)}_h$ be the $i$-th component of the state vector $\mathbf{U}_h$ and $N$ be the total number of components. 
The coefficient $\sigma^m_{e, K} \left(\mathbf{U}_h\right)$ is defined as
\begin{equation}\label{sigma}
	\sigma^m_{e, K} \left(\mathbf{U}_h\right)= \max_{1 \leq i \leq N} \sigma^m_{e, K} \left(u^{\left(i\right)}_h\right),
\end{equation}
with 
\begin{equation} \label{sigma-e}
	\sigma^m_{e, K} \left(u^{\left(i\right)}_h\right) = \begin{dcases}
		0,  & \mathrm{if} \ u^{\left(i\right)}_h \equiv \operatorname{avg}\left(u_h^{\left(i\right)}\right),
		\\
		\dfrac{\left(2m+1\right)h^m_{e, K}}{2\left(2k-1\right)m!}  \dfrac{\sum\limits_{\left| \boldsymbol{\alpha} \right|=m}\dfrac{1}{\left|e\right|} \int_e  \big| [\![ \partial^{\boldsymbol{\alpha}} u_h^{(i)} ]\!] \big|  ds}{\|u_h^{(i)}-\operatorname{avg}(u_h^{(i)})\|_{L^{\infty}(\Omega)}}, & \mathrm{otherwise},
	\end{dcases}
\end{equation}
where $\partial^{\boldsymbol{\alpha}} u_h^{(i)}$ is the partial derivative defined as 
\begin{equation*}
	\partial^{\boldsymbol{\alpha}}u_h^{(i)}(\boldsymbol x) = \frac{\partial ^{|\boldsymbol{\alpha}|}}{\partial x^{\alpha_1}\partial y^{\alpha_2}}u_h^{(i)}(\boldsymbol x),
\end{equation*}
$\big| [\![ \partial^{\boldsymbol{\alpha}} u_h^{(i)} ]\!]\big| $ is the absolute value of the jump of $\partial^{\boldsymbol{\alpha}} u_h^{(i)}$, and $\operatorname{avg}\left(u^{\left(i\right)}_h\right)= \dfrac{1}{\left|\Omega\right|}\int_{\Omega} u^{\left(i\right)}_h \ d\boldsymbol{x}$ is the average of $u^{\left(i\right)}_h$ over $\Omega$. The common scalar quantity $\|u_h^{(i)}-\operatorname{avg}(u_h^{(i)})\|_{L^{\infty}(\Omega)}$ is computed at the beginning of the OE procedure and then used for all the elements. The face integral of $\big| [\![ \partial^{\boldsymbol{\alpha}} u_h^{(i)} ]\!] \big| $ on $e$ is calculated by using Gauss quadrature rules.

Computation of the damping coefficients on a rectangular element $K= \left[ x_{1}, x_{2} \right] \times \left[ y_{1}, y_{2} \right] $ is presented as an illustration of how to use formulae \eqref{delta} to \eqref{sigma-e}. Let $\beta_{K}^x$, $\beta_{K}^y$ be the estimates of the local maximum wave speeds in the $x$- and $y$-directions, respectively. Then $\delta_K^m\left(\mathbf{U}_h\right)$ can be calculated by
\begin{equation*}
	\begin{aligned}
	&\delta_K^m\left(\mathbf{U}_h\right)\\
	&=\max _{1 \leq q \leq N}\left(\frac{\beta_{K}^x\left(\sigma_{x=x_2}^m(u_h^{(q)})+\sigma_{x=x_1}^m(u_h^{(q)})\right)}{h_{x}}+\frac{\beta_{K}^y\left(\sigma_{y=y_2}^m(u_h^{(q)})+\sigma_{y=y_1}^m(u_h^{(q)})\right)}{h_{y}}\right),
	\end{aligned}
\end{equation*}
with
\begin{equation*}
	\begin{aligned}
		& \sigma_{x=x_2}^m(u_h^{(q)})= \begin{cases}0, & \text { if } u_h^{(q)} \equiv \operatorname{avg}(u_h^{(q)}), \\
			\dfrac{ \frac{2m+1}{2(2k-1)m!} h_{x}^m\sum\limits_{|\boldsymbol \alpha|=m} \frac{1}{h_{y}}  \int_{y_1}^{y_2}\big|[\![ \partial^{\boldsymbol{\alpha}} u_h^{(q)} ]\!]\big|_{x=x_2} \mathrm{d} y}{{\|u_h^{(q)}-\operatorname{avg}(u_h^{(q)})\|_{L^{\infty}(\Omega)}}},& \text { otherwise, }\end{cases} \\
		& \sigma_{y=y_2}^m(u_h^{(q)})= \begin{cases}0, & \text { if } u_h^{(q)} \equiv \operatorname{avg}(u_h^{(q)}), \\
			\dfrac{ \frac{2m+1}{2(2k-1)m!} h_{y}^m\sum\limits_{|\boldsymbol \alpha|=m} \frac{1}{h_{x}} \int_{x_1}^{x_2}\big|[\![ \partial^{\boldsymbol \alpha} u_h^{(q)} ]\!]\big|_{y=y_2} \mathrm{d} x}{{\|u_h^{(q)}-\operatorname{avg}(u_h^{(q)})\|_{L^{\infty}(\Omega)}}}, & \text { otherwise. }\end{cases}
	\end{aligned}
\end{equation*}

It is observed from \eqref{damping-problem} and \eqref{delta} - \eqref{sigma-e} that, the damping term $-\Sigma \left(\mathbf{U}_h\right) \mathbf{U}_\sigma$ is a high-order term in smooth areas, thus the seemingly \textquote{first-order} operator splitting presented in Section \ref{subsec:OEDG} does not affect the order of accuracy of the OEDG scheme. Actually, it is proved in \cite{Peng2023OEDG} that, the OEDG method can preserve order of accuracy in smooth regions. The accuracy-preserving property of the OEDG method is also demonstrated by the numerical results in Section \ref{Sect4} of this paper. The strong discontinuity capturing capability of the OEDG method is demonstrated by the numerical results in Section \ref{Sect4}. It is observed in \eqref{delta} that the damping coefficients are proportional to wave speeds, leading to an evolution-invariant property of the OEDG method. The OEDG method is also scale-invariant, as the damping coefficients are normalized by the $L^\infty$-norm of the solution fluctuation. 

\subsection{Locally divergence-free oscillation-eliminating discontinuous Galerkin method}

There are two main difficulties in developing high-order DG schemes for ideal compressible MHD equations, namely the preservation of magnetic divergence-free property and the elimination of spurious oscillations near discontinuities. The LDF-DG scheme presented in Section \ref{Sect2} can be used to obtain high-order accurate divergence-free numerical solutions of ideal compressible MHD equations. The OE procedure presented in Section \ref{Sect3} can be used to suppress spurious oscillations in the vicinity of discontinuities while preserve accuracy in smooth regions. Therefore, we propose a LDF-OEDG method to obtain high-order accurate divergence- and oscillation-free numerical solutions of ideal compressible MHD equations. The implementation of the LDF-OEDG scheme is described in Algorithm \ref{alg:LDF-OEDG}. 

\begin{algorithm}[htbp!]
	\caption{LDF-OEDG for ideal compressible MHD equations.}  
	\label{alg:LDF-OEDG}
	\begin{algorithmic}[1] 

		\Function {$\mathbf{U}^{n+1}_\sigma$=LDF-OEDG}{$\mathbf{U}^{n}_\sigma$, $\tau$} 
		
		\State Set $\mathbf{U}^{n, 0}_\sigma= \mathbf{U}^{n}_\sigma$
		
		\State Set SSPRK3 coefficients: $c_1=1$, $c_1=\frac{1}{4}$, $c_3=\frac{2}{3}$
		
		\For{$s \gets 1, 3$}
		\State Surface and volume flux integrals in \eqref{semi-DG} to compute right-hand-side $\mathcal{T}_f \left(\mathbf{U}^{n, s-1}_\sigma\right)$ 
		
		\State Runge-Kutta stage update $\mathbf U_h^{n,s} = \left(1-c_s\right) \mathbf U_\sigma^n + c_s \left(\mathbf U_\sigma^{n,s-1}+\tau \mathcal T_f\left(\mathbf U_\sigma^{n,s-1}\right)\right)$
		
		\State OE procedure $ \mathbf U_\sigma^{n,s} =\mathcal{F}_\tau \mathbf U_h^{n,s}$ using the exact damping operator in \eqref{OE-solution}
		
		\State Splitting of $\mathbf U_\sigma^{n,s}$ into two parts: the flow field $\mathbf{U}^{n,s}_{\mathrm{\sigma, F}}$ and the magnetic filed $\mathbf{B}^{n,s,*}_\sigma$
		
		\State Projection of $\mathbf{B}^{n,s,*}_\sigma$ to obtain a LDF magnetic filed $\mathbf{B}^{n,s}_\sigma$, as in \eqref{formula-projection}
		
		\State Formation of a divergence- and oscillation-free solution $ \mathbf U_\sigma^{n,s} = \left(\mathbf{U}^{n,s}_{\mathrm{\sigma, F}}, \mathbf{B}^{n,s}_\sigma\right)^{\mathrm{T}}$
		\EndFor
		
		\State Update solution $\mathbf{U}^{n+1}_\sigma= \mathbf{U}^{n,3}_\sigma$

		\EndFunction
		
	\end{algorithmic}  
\end{algorithm} 

\section{Numerical results}\label{Sect4}
This section presents numerical results for several benchmark two-dimensional ideal MHD cases to verify the accuracy, resolution and robustness of the proposed LDF-OEDG method. We specifically explore the third-order accurate ($k=2$) LDF-OEDG method on rectangular meshes. The third-order SSP Runge-Kutta scheme is used in time integration, with the time step size
\begin{equation*}
	\tau= \dfrac{\mathrm{CFL}}{\lambda_x/h_x + \lambda_y/h_y},
\end{equation*}
where $\lambda_x$ and $\lambda_y$ are the largest wave speeds in $x$- and $y$-direction, respectively. The CFL number is taken as 0.15. 
For all cases, the adiabatic index is $\gamma = 5/3$, unless otherwise noted.

\begin{example}\label{exam1}
	(Accuracy test) The first example is a two-dimensional vortex problem, which was adapted to MHD equations by Balsara \cite{Balsara2004Second} from the isentropic vortex problem \cite{Hu1999weighted} in gas dynamics.  The solution is a smooth vortex stably convected with the velocity and magnetic field. The mean flow is $(\rho,\mathbf{u},\mathbf{B},p)=(1,1,1,0,0,1)$. A vortex is added to the mean flow with the following perturbations
	\begin{equation*}
		\begin{aligned}
			&(\delta u_x, \delta u_y) = \frac{1}{2\pi} e^{0.5(1-r^2)} \left(-y,x\right),\\
			&(\delta B_x, \delta B_y) = \frac{1}{2\pi}e^{0.5(1-r^2)}\left(-y,x\right),\\
			&\delta p = \frac{-r^2}{8\pi^2}e^{(1-r^2)},
		\end{aligned}
	\end{equation*}
	where $r^2=x^2+y^2$. The computational domain is $[-5,5] \times [-5,5]$. Periodic boundary conditions are imposed on the left/right and top/bottom boundaries, respectively. 
\end{example}
This example is used to assess the order of convergence of the LDF-OEDG method for smooth problems. The simulations are performed up to $t=20$, on a set of successively refined rectangular meshes.
The $L^2$-errors and orders of accuracy for some representative conservative variables are listed in Table \ref{ErrTab1}. The test results show that, the proposed LDF-OEDG method with $k=2$ can achieve the expected third-order accuracy.

\begin{table}[hbtp!]
	\caption{Accuracy test results for the vortex problem.}
	\centering
	\begin{tabular}{ccccccccccc}
		\toprule
		\centering
		& \multicolumn{2}{c}{$\rho$} & \multicolumn{2}{c}{$\rho u_x$} & \multicolumn{2}{c}{$B_x$} & \multicolumn{2}{c}{$E$}\\
		\midrule
		Mesh & $L^2$-error & Order & $L^2$-error & Order & $L^2$-error & Order & $L^2$-error & Order\\
		\midrule
		$16\times 16$ & 2.12E-03  & & 1.10E-02 &  & 1.14E-02 &  & 1.50E-02 &\\
		$32\times 32$ & 2.59E-04 & 3.03 & 5.92E-04 & 4.22 & 5.84E-04 & 4.28 & 9.88E-04 &3.93\\
		$64\times 64$ & 4.45E-05 & 2.54& 5.54E-05 & 3.42 & 5.17E-05 & 3.50 & 1.23E-04 &3.01\\
		$128\times 128$ & 7.29E-06 & 2.61 & 7.83E-06 & 2.82 & 6.51E-06  & 2.99 & 1.86E-05 &2.73\\
		\bottomrule
	\end{tabular}\label{ErrTab1}
\end{table}

\begin{figure}[hbtp!]
	\centering
	\begin{subfigure}{0.48\textwidth}
		\includegraphics[width=\linewidth]{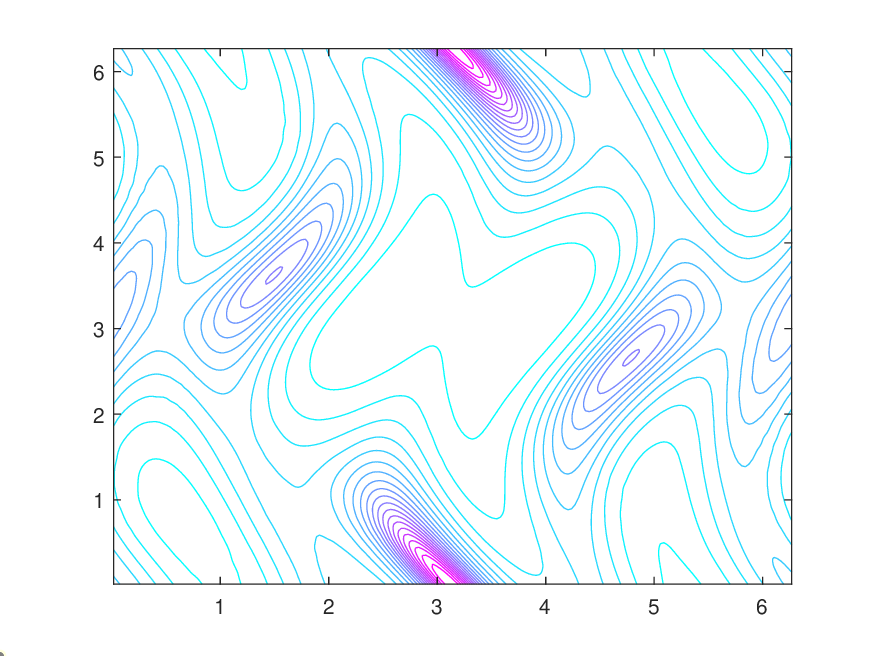}
		\caption{$t=0.5$}
		\label{fig:OEOTV0_5}
	\end{subfigure}
	\begin{subfigure}{0.48\textwidth}
		\includegraphics[width=\linewidth]{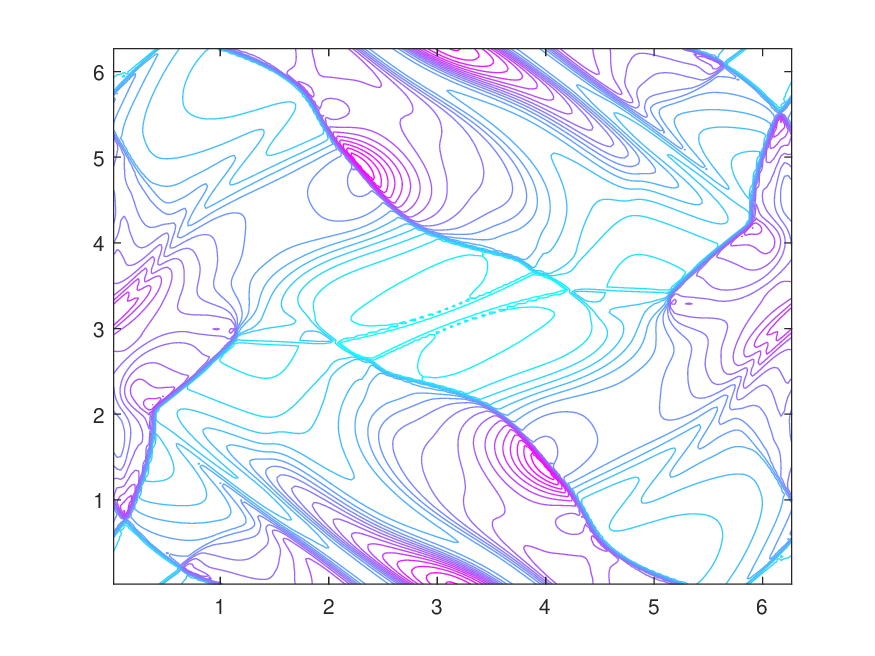}
		\caption{$t=2$}
		\label{fig:OEOTV2}
	\end{subfigure}
	\begin{subfigure}{0.48\textwidth}
		\includegraphics[width=\linewidth]{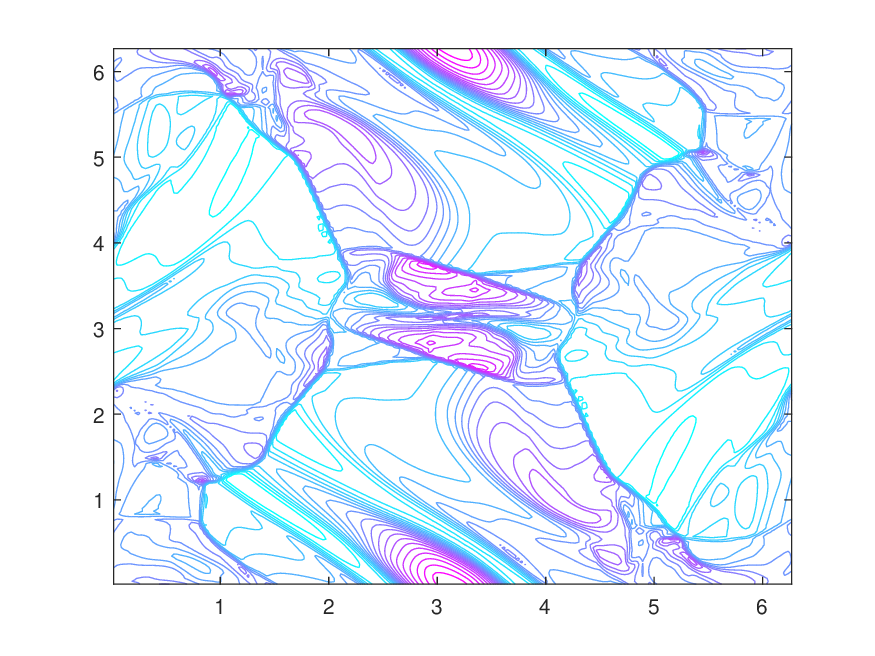}
		\caption{$t=3$}
		\label{fig:OEOTV3}
	\end{subfigure}
	\begin{subfigure}{0.48\textwidth}
		\includegraphics[width=\linewidth]{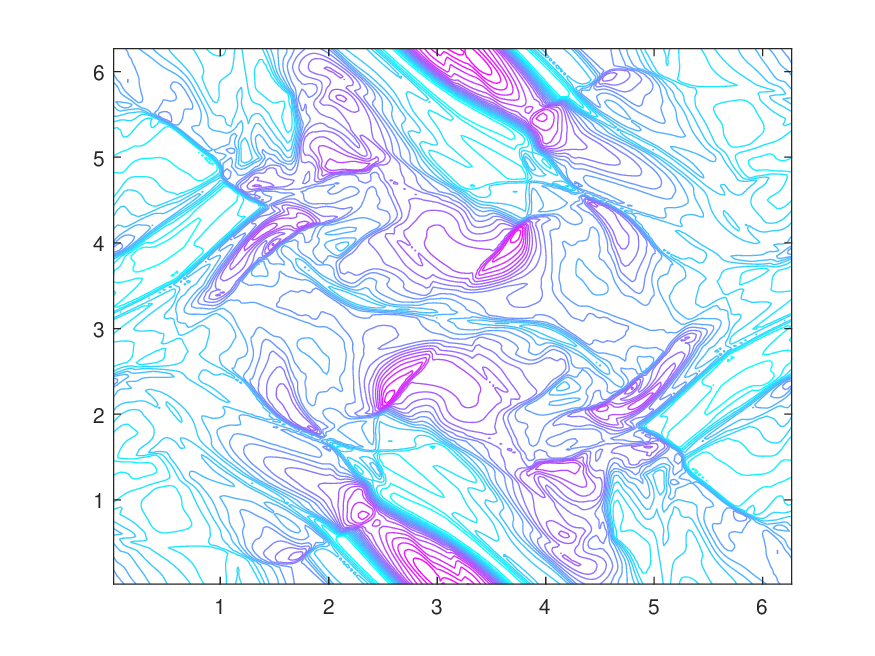}
		\caption{$t=4$}
		\label{fig:OEOTV4}
	\end{subfigure}
	\caption{Orszag-Tang problem. Density contours computed by LDF-OEDG on a $192\times 192$ mesh.}
	\label{OEOTVPlot}
\end{figure}

\begin{example}\label{exam2}
	(Orszag-Tang problem)
	The classic Orszag-Tang problem \cite{Orszag1979Small-scale} is widely used to assess resolution of numerical methods. The initial conditions are
	$$
	(\rho,\mathbf{u},\mathbf{B},p)=(\gamma^2,-\sin y,\sin x,-\sin y,\sin 2x, \gamma).
	$$
	The computational domain is $[0,2\pi]\times[0,2\pi]$, with periodic boundary conditions imposed on top/bottom and left/right boundaries. The computational domain is divided into $192\times 192$  uniform cells in this test. 
\end{example}
Figure \ref{OEOTVPlot} plots the density contours at $t=0.5,$ $2,$ $3$ and $4$ computed by the LDF-OEDG method. One can observe that the solution is smooth at early stage ($t=0.5$). Then the shocks have already appeared at $t=2$. After that, the shocks interact with each other and the complicated structures involving multiple shocks emerge at $t=3$ and $4$. It is observed that the results are consistent with those in \cite{Jiang1999A,Li2005Locally},  and the complex shock structures are correctly captured with high resolution by the proposed LDF-OEDG method.

\begin{example}\label{exam3}
	(Rotor problem)
	The rotor example \cite{Balsara1999A,Toth2000The} is used to further verify the resolution of the LDF-OEDG method. In this problem, a dense rotating disk of fluid is located in the central area, surrounded by the ambient fluid. An intermediate flow state is set between these two areas for a smooth initial condition, which is given by
	$$
	(\rho,\mathbf{u},\mathbf{B},p)=(\rho^0,u^0_x,u^0_y,2.5/\sqrt{4\pi},0,0.5)
	$$
	with
	\begin{align*}
		\begin{split}
			(\rho^0,u^0_x,u^0_y)=
			\begin{cases}
				(10,~-(y-0.5)/r_0,~(x-0.5)/r_0),\ &if~r \le r_0,\\
				(1+9f,~-f\cdot(y-0.5)/r,~f\cdot(x-0.5)/r),\ &if~r_0<r<r_1,\\
				(1,~0,~0),\ &if~r \ge r_1,
			\end{cases}
		\end{split}
	\end{align*}
	where $r_0=0.1,$ $r_1=0.115,$ $f=(r_1-r)/(r_1-r_0)$ and $r=\left[(x-0.5)^2+(y-0.5)^2\right]^{1/2}$. The computational domain is $[0,1]\times[0,1]$, with periodic boundary conditions imposed on top/bottom and left/right boundaries. 
\end{example}

\begin{figure}[htbp!]
	\centering
	\includegraphics[scale=0.35]{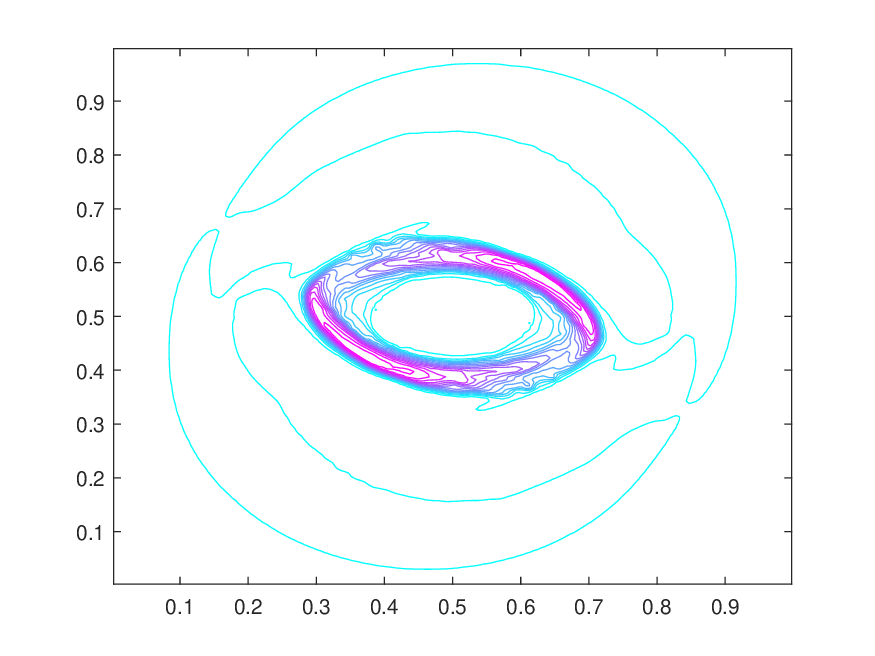}
	\includegraphics[scale=0.35]{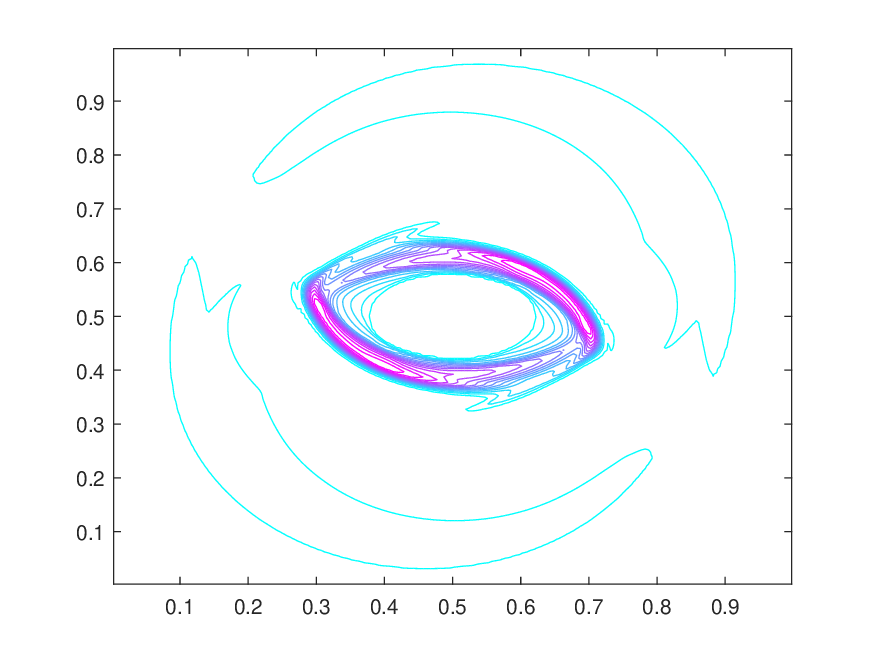}\\
	\includegraphics[scale=0.35]{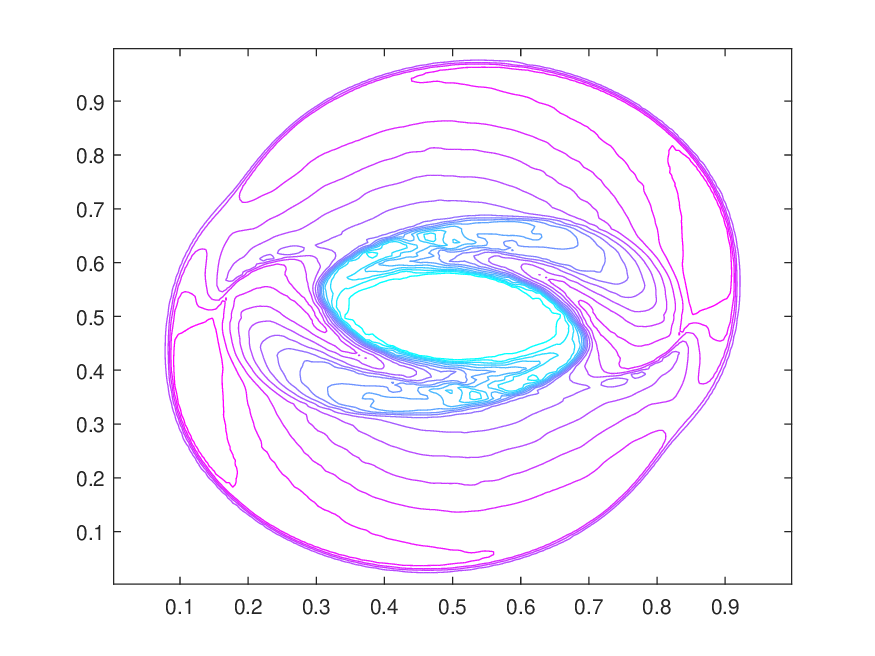}
	\includegraphics[scale=0.35]{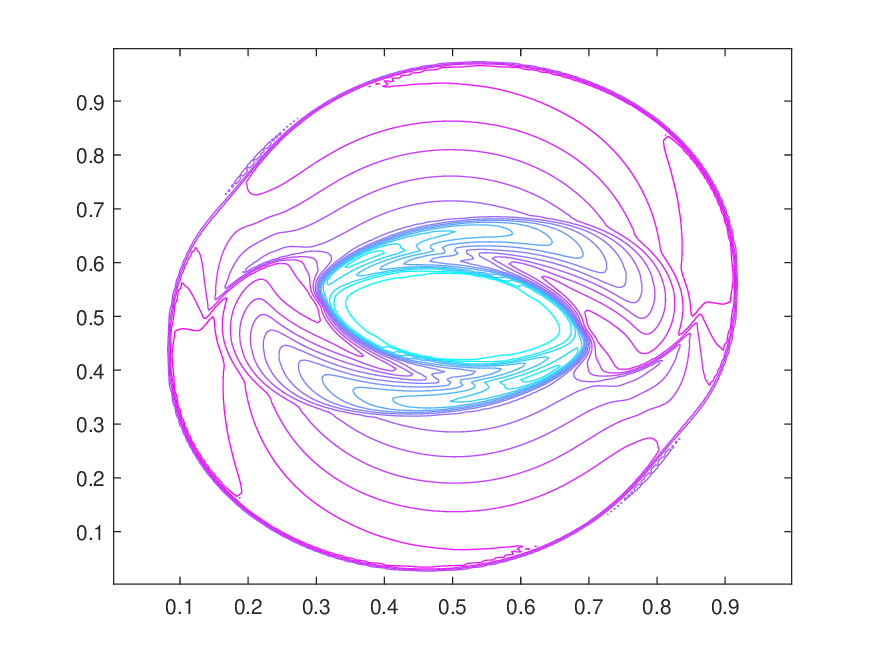}\\
	\includegraphics[scale=0.35]{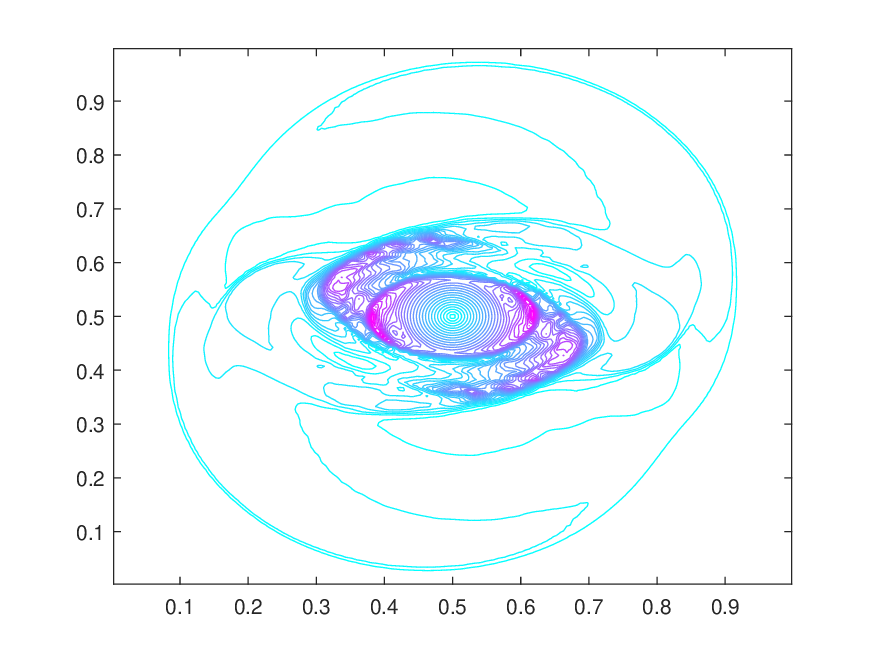}
	\includegraphics[scale=0.35]{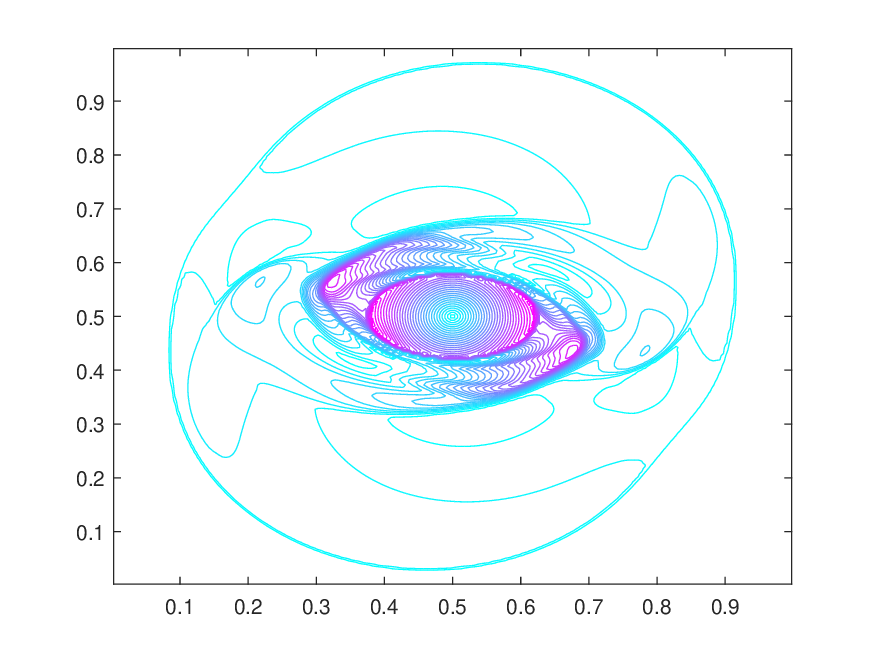}\\
	\includegraphics[scale=0.35]{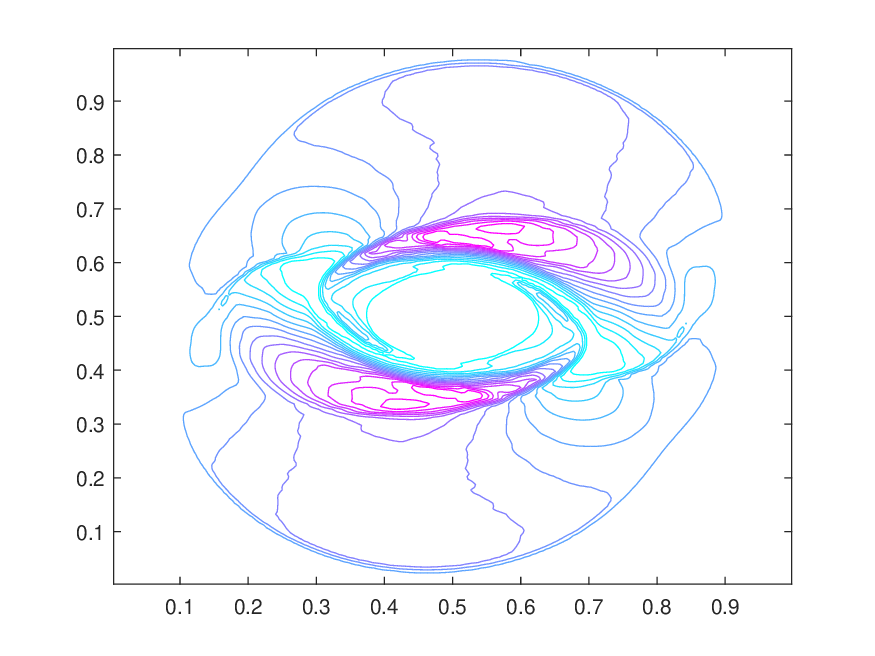}
	\includegraphics[scale=0.35]{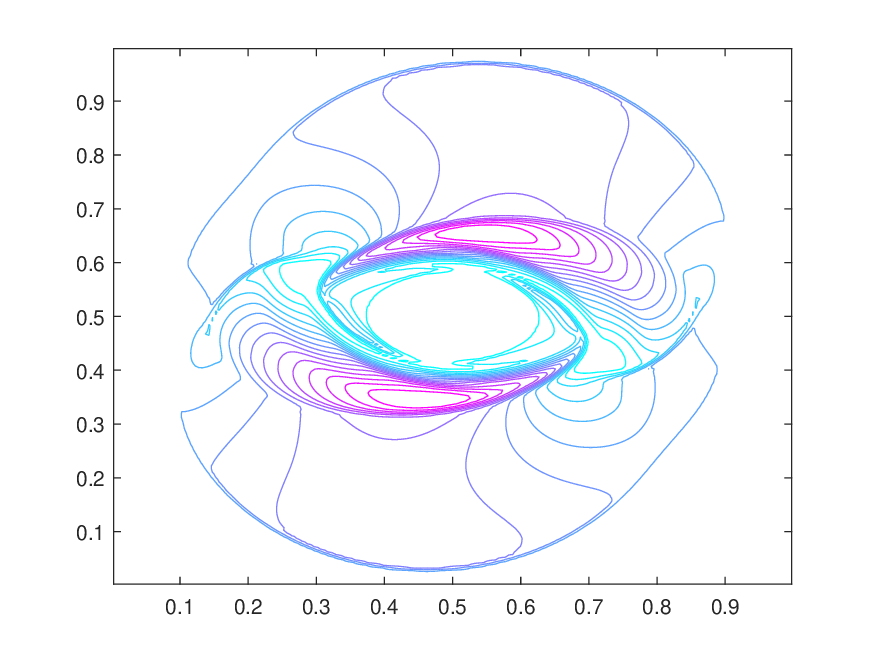}\\
	\caption{Rotor problem. From top to bottom: contour plots of the density $\rho$,  thermal pressure $p$, Mach number $\left\|\mathbf{u}\right\|_2/c$ and magnetic pressure $\left\|\mathbf{B}\right\|_2^2/2$, respectively. The solutions are computed on a $200 \times 200$ mesh. Left: LDF-DG with TVB limiter; right: LDF-OEDG.}\label{RotorPlot1}
\end{figure}

Numerical results computed by the LDF-OEDG method at $t=0.295$ on a $200\times 200$ mesh, including the contours of density $\rho$, thermal pressure $p$, hydrodynamic Mach number $\left\|\mathbf{u}\right\|_2/c$ with $c=\sqrt{\frac{\gamma p}{\rho}}$, and magnetic pressure $\left\|\mathbf{B}\right\|_2^2/2$, are shown in Figure \ref{RotorPlot1}. Results computed by a LDF-DG method using the TVB limiter is also shown in Figure \ref{RotorPlot1} for a comparison. It is observed in Figure \ref{RotorPlot1} that, the flow fields in the central area computed by the LDF-DG method with the TVB limiter are more noisy than those computed by the proposed LDF-OEDG method, demonstrating the superior shock capturing capability of the LDF-OEDG method.

\begin{figure}[htbp!]
	\centering
	\includegraphics[scale=0.35]{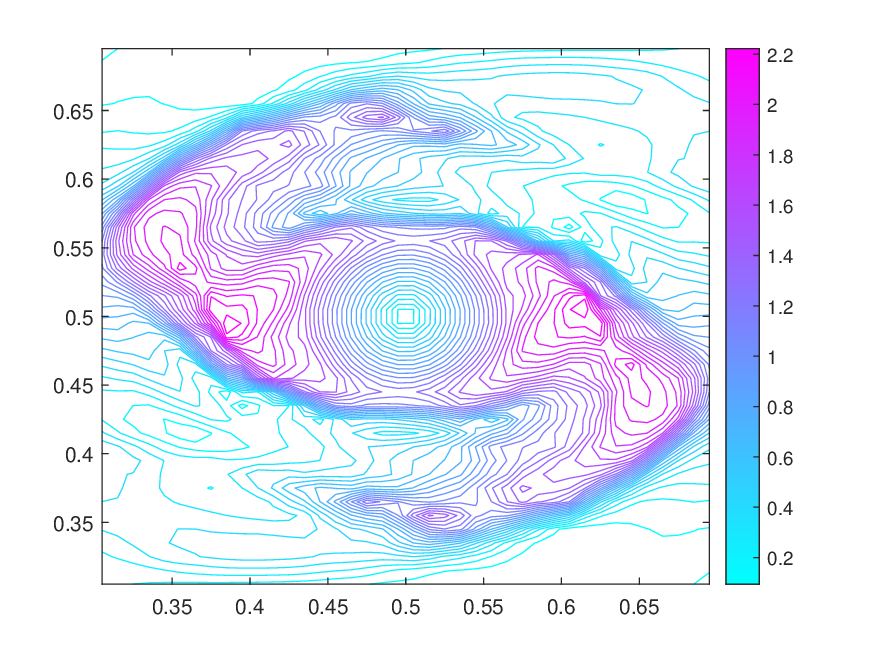}
	\includegraphics[scale=0.35]{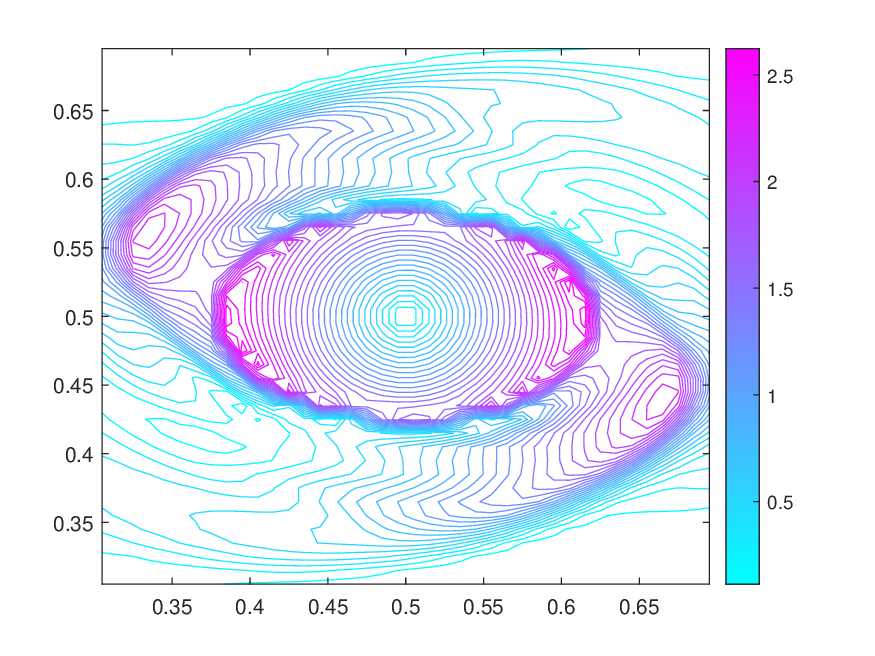}\\
	\includegraphics[scale=0.35]{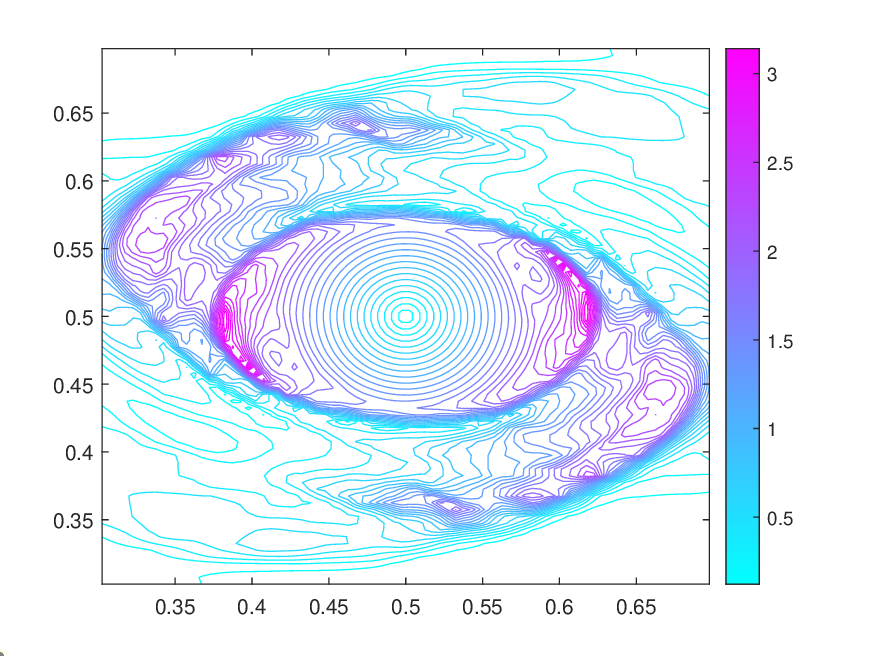}
	\includegraphics[scale=0.35]{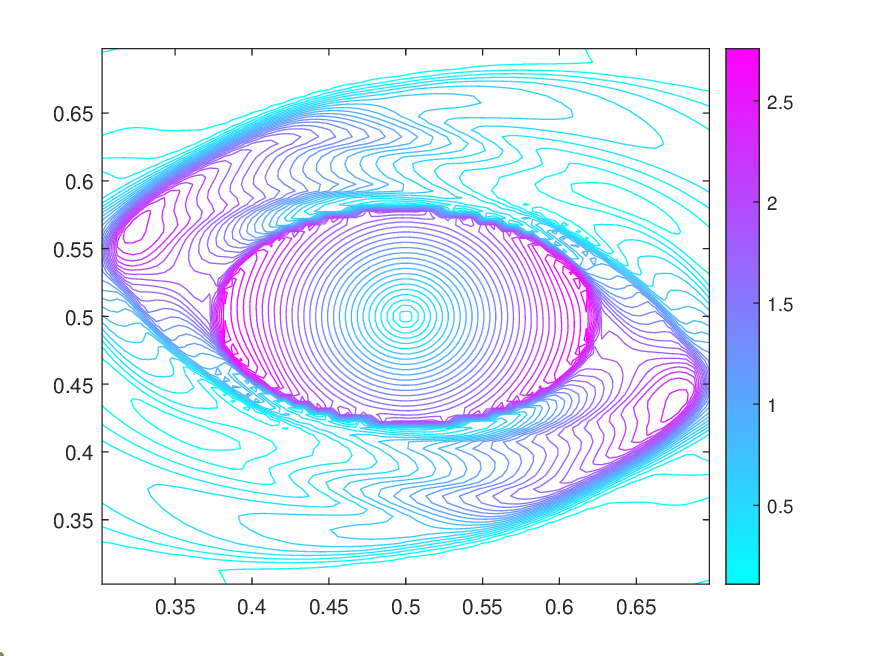}\\
	\includegraphics[scale=0.35]{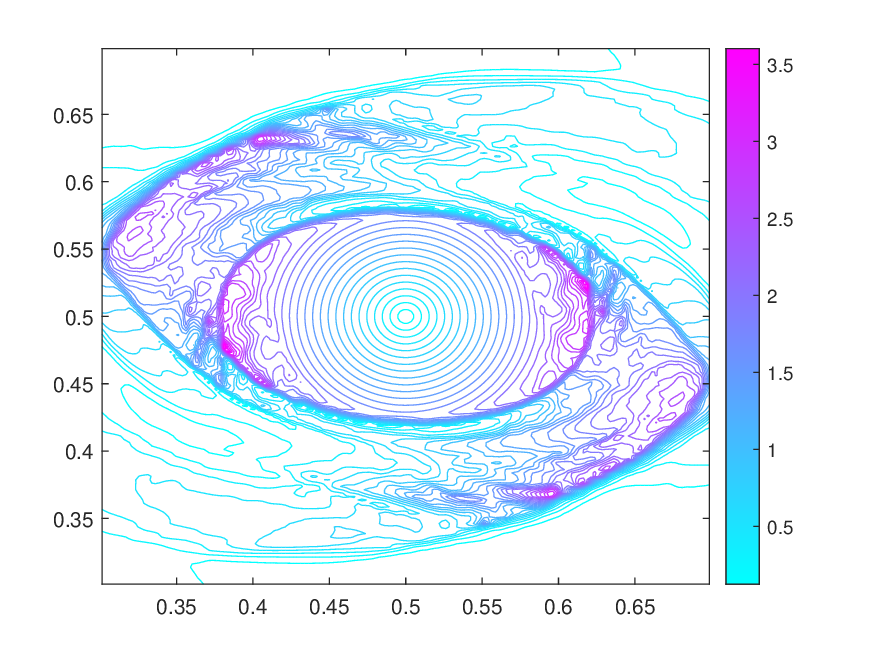}
	\includegraphics[scale=0.35]{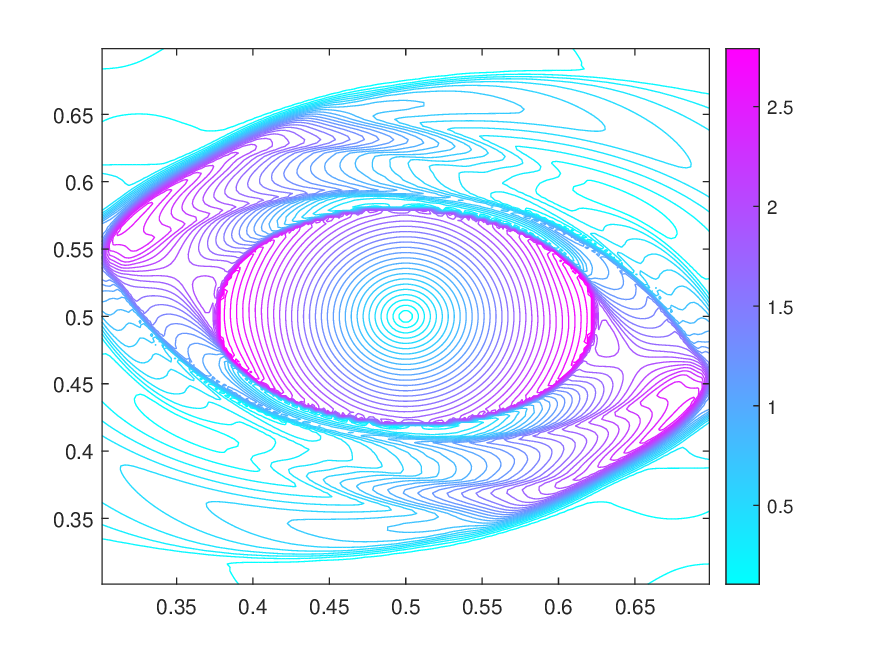}\\
	\includegraphics[scale=0.35]{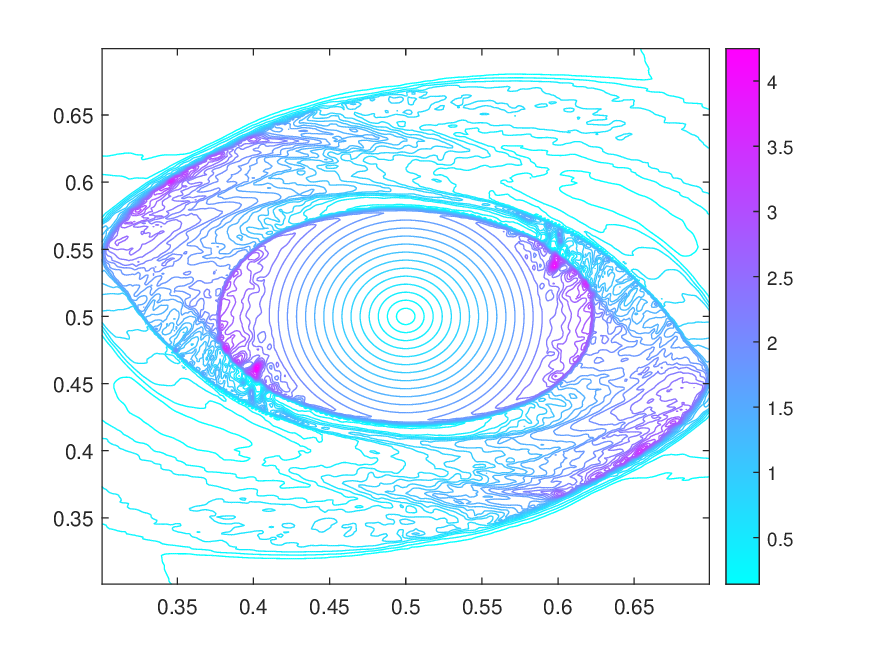}
	\includegraphics[scale=0.35]{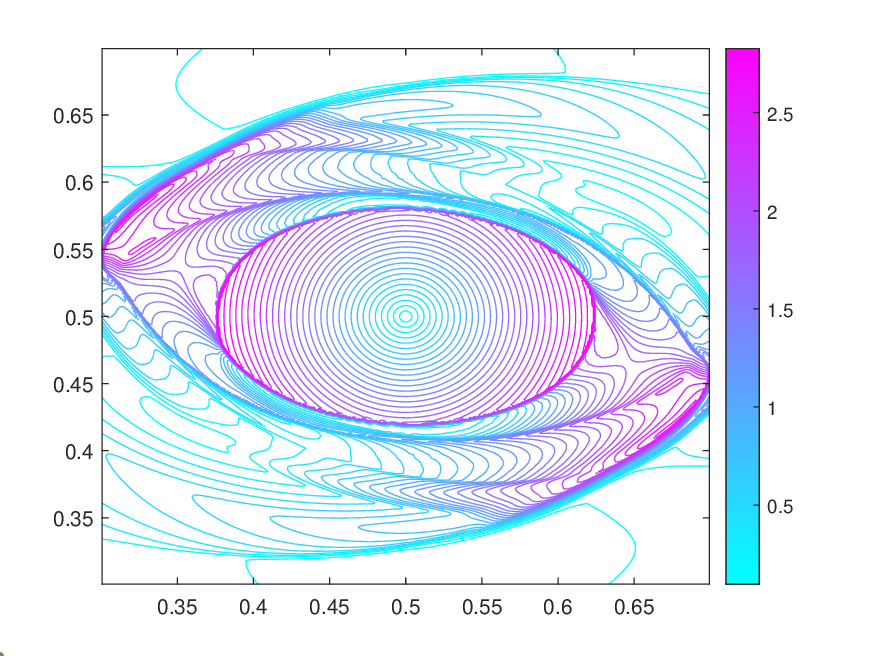}\\
	\caption{Rotor problem. Mach number contours (zoom-in central part). From top to bottom: on meshes with $100\times 100$, $200\times 200$, $400\times 400$ and $800 \times 800$ cells, respectively. Left: LDF-DG with TVB limiter; right: LDF-OEDG.}\label{RotorPlot}
\end{figure}

For a more detailed resolution comparison, in Figure \ref{RotorPlot}, we present the Mach number contours in the central rotating area computed on four successively refined rectangular meshes with $100\times 100$, $200\times 200$, $400\times 400$ and $800\times 800$ cells, respectively. It is observed from Figure \ref{RotorPlot} that, the solution resolution increases with the grid refinement. Furthermore, the solutions computed by the LDF-OEDG method have remarkably higher resolution than those computed by the LDF-DG method with the TVB limiter, indicating the superior performance of the OE procedure to the TVB limiter.

\begin{example}\label{exam4}
	(Blast problem)
	The classic MHD blast wave problem \cite{Balsara1999A} is used to assess the shock capturing capability and robustness of the LDF-OEDG method. This a challenging problem, in which a strong circular magneto-sonic shock formulates and propagates into the ambient plasma with a small plasma-beta $\beta=2p/{\left\|\mathbf{B}\right\|_2^2}$. A smaller $\beta$ yields a higher probability of producing negative pressure. Initially,
	\begin{align*}
		\begin{split}
			(\rho,\mathbf{u},\mathbf{B},p)=
			\begin{cases}
				\left(1,0,0,\frac{100}{\sqrt{4\pi}},0,1000\right),\ &if~r\leq R,\\
				\left(1,0,0,\frac{100}{\sqrt{4\pi}},0,0.1\right),\ &if~r>R,\\
			\end{cases}
		\end{split}
	\end{align*}
	where $R=0.1,$ $r=\left(x^2+y^2\right)^{1/2}$. This initial condition corresponds to a very small plasma-beta $\beta=2.51E-04$, making the numerical simulation quite challenging. The adiabatic index is $\gamma=1.4$. The computational domain is $[-0.5, 0.5]\times[-0.5, 0.5]$. Outflow boundary conditions are imposed on domain boundaries. 
\end{example}

\begin{figure}[htbp]
	\centering
	\begin{subfigure}{0.48\textwidth}
		\centering
		\includegraphics[width=\linewidth]{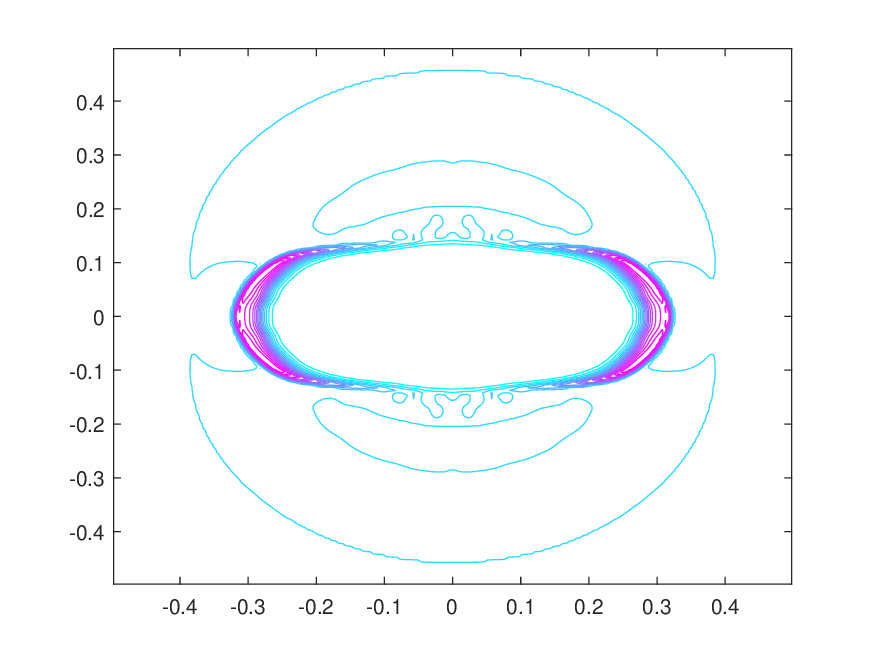}
		\caption{$\rho$}
	\end{subfigure}%
	\begin{subfigure}{0.48\textwidth}
		\centering
		\includegraphics[width=\linewidth]{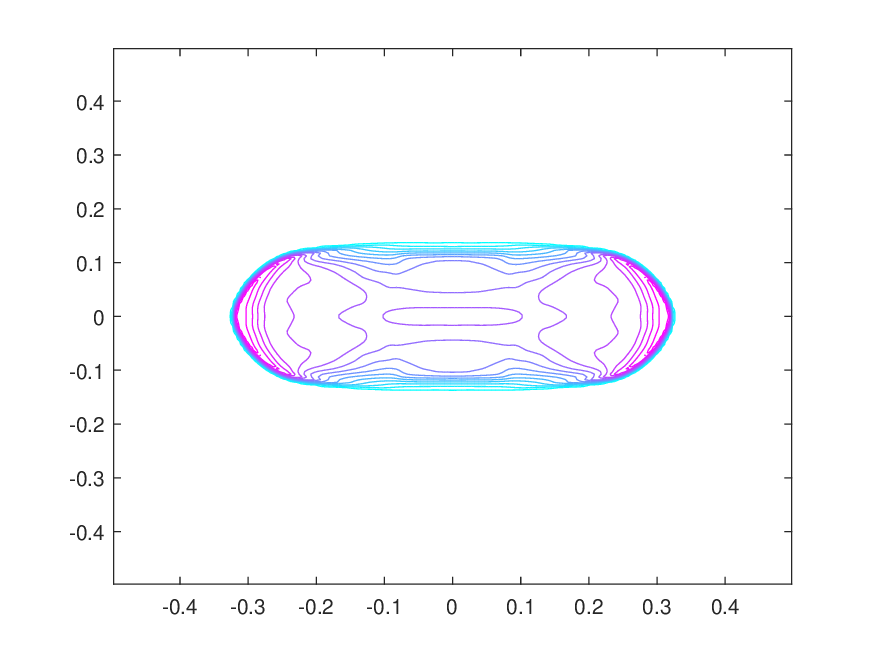}
		\caption{$p$}
	\end{subfigure}
	\begin{subfigure}{0.48\textwidth}
		\centering
		\includegraphics[width=\linewidth]{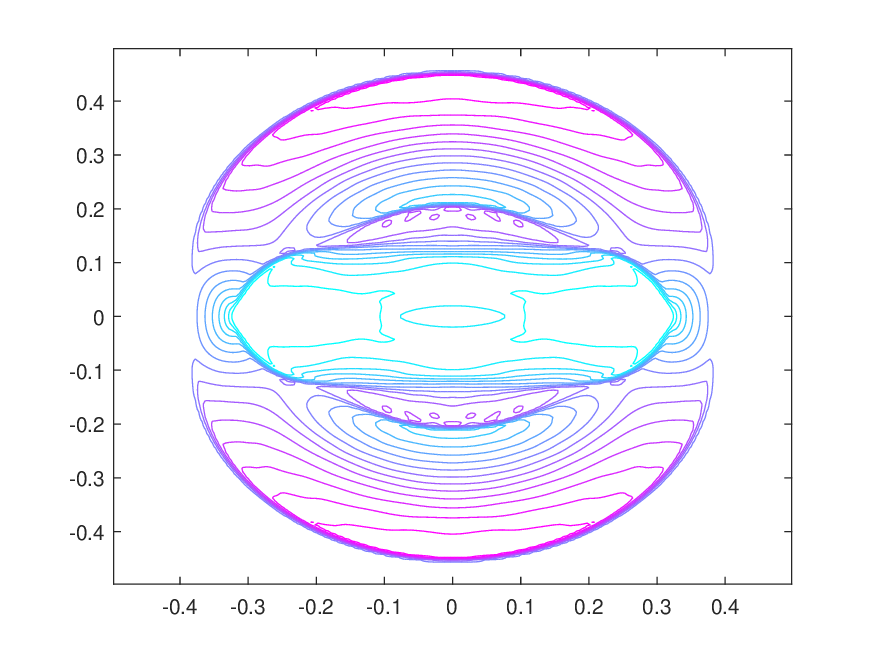}
		\caption{$\left\| \mathbf{B}\right\|_2^2/2$}
	\end{subfigure}
	\caption{Blast problem. Contours of density, thermal pressure, and magnetic pressure computed by LDF-OEDG.}\label{BlastPlot}
\end{figure}

A simulation is performed on a $200 \times 200$ mesh until $t= 0.01$ using the LDF-OEDG method. As the plasma-beta $\beta$ is quite small, the numerical scheme needs to be robust, otherwise there is a high probability of simulation failure. The numerical results, including the contours of density $\rho$, thermal pressure $p$ and magnetic pressure $\left\| \mathbf{B}\right\|_2^2/2$, are shown in Figure \ref{BlastPlot}. It is observed in Figure \ref{BlastPlot} that the simulation results are in good agreement with those reported in \cite{Balsara1999A,Christlieb2015Positivity}, demonstrating the strong shock capturing capability and robustness of the LDF-OEDG method.

\begin{example}\label{exam5}
	(Loop advection problem) The magnetic field loop advection problem, originally proposed in \cite{Gardiner2005An}, is used to assess the capability of a numerical scheme to preserve the loop's initial shape and magnetic field strength over time.
	The initial condition is
	\begin{equation*}
		(\rho,\mathbf{u},\mathbf{B},p)=(1,2,1,B_x,B_y,1),
	\end{equation*}
	where
	\begin{equation*}
		\left(B_x,B_y\right)= \left(\dfrac{\partial A_z}{\partial y}, -\dfrac{\partial A_z}{\partial x}\right), \quad 	
		\begin{aligned}
			\begin{split}
				A_z=
				\begin{cases}
					A_0(R-r),\ &if~r\leq R,\\
					0,\ &if~r>R,\\
				\end{cases}
			\end{split}
		\end{aligned}
	\end{equation*}
	with $A_0=10^{-3}$, $R=0.3$ and $r=\left(x^2+y^2\right)^{1/2}$. It is noted that $B_x$ and $B_y$ are the first two components of the curl of the magnetic vector potential $\nabla\times(0,0,A_z)^{\mathrm{T}}$.
	The computational domain is $[-1,1]\times[-0.5,0.5]$, with periodic boundary conditions imposed on top/bottom and left/right boundaries. 
\end{example}

\begin{figure}[htbp]
	\centering
	\begin{subfigure}{0.48\textwidth}
		\centering
		\includegraphics[width=\linewidth]{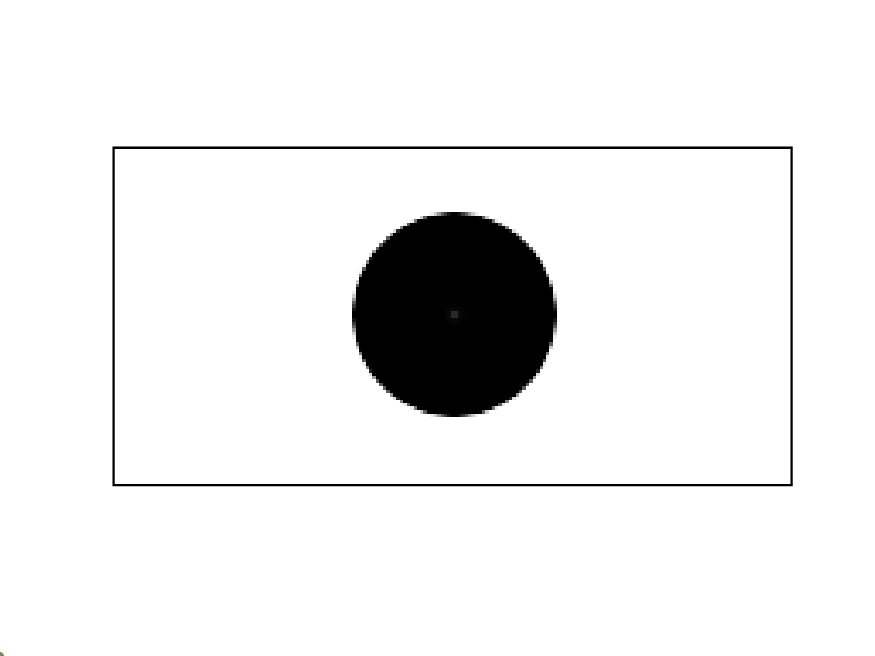}
		\caption{$t=0$}
	\end{subfigure}%
	\begin{subfigure}{0.48\textwidth}
		\centering
		\includegraphics[width=\linewidth]{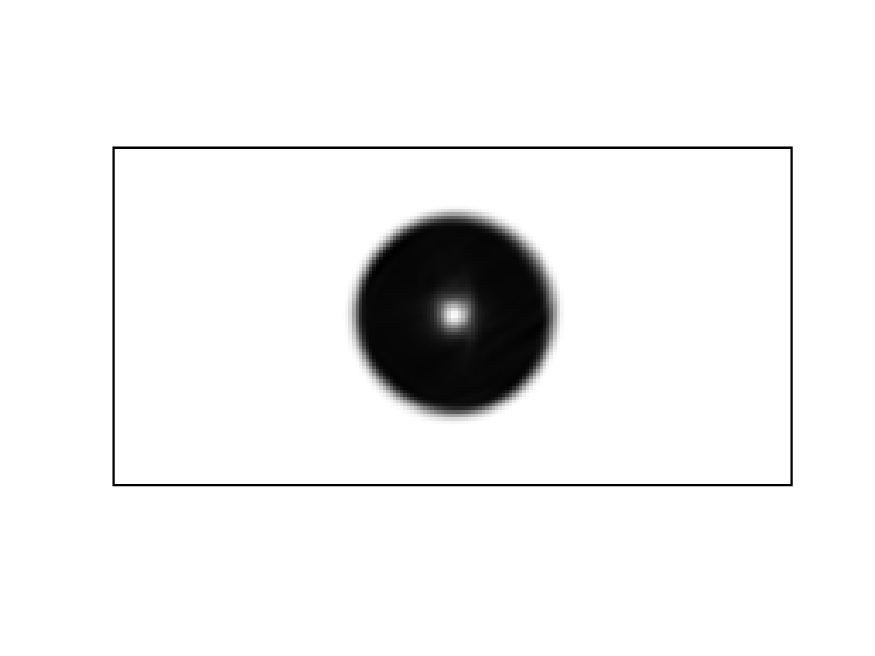}
		\caption{$t=2$}
	\end{subfigure}
	\begin{subfigure}{0.48\textwidth}
		\centering
		\includegraphics[width=\linewidth]{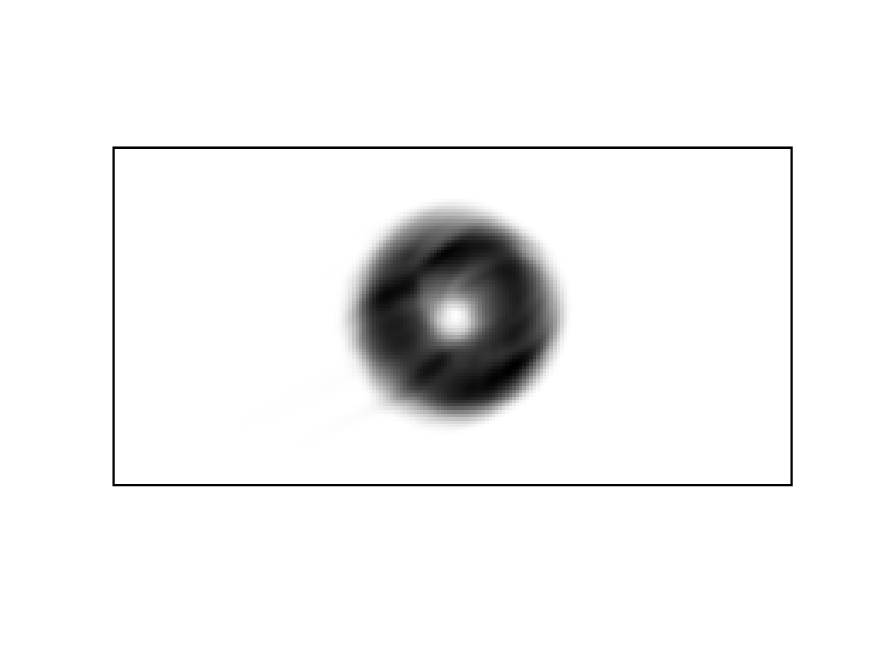}
		\caption{$t=10$}
	\end{subfigure}
	\caption{Loop advection problem. Gray-scale images of $\left\| \mathbf{B} \right\|_2^2$.}
	\label{LoopPlot} 
\end{figure}

The gray-scale images of $\left\| \mathbf{B} \right\|_2^2$ at $t = 0$, $2$ and $10$, computed by the LDF-OEDG method on a $200 \times 100$ mesh, are shown in Figure \ref{LoopPlot}. It is observed that the majority of field dissipation occurs at the center and boundaries of the field loop, and the loop's initial shape is preserved over time. These phenomena agree with the results in \cite{Gardiner2005An,Fu2018Globally}. 

\begin{example}\label{exam6}
	(Shock cloud interaction problem)
	The last example is a cloud–shock interaction \cite{Dai1998A}, which is a model widely used in astrophysics. Its main characteristic is the interaction between a high density cloud and a strong shock wave.
	Set
	\begin{align*}
		\begin{split}
			\Omega_1=&\{(x,y)|0\leq x \leq 1.2,~0\leq y\leq 1\},\\
			\Omega_2=&\{(x,y)|1.2\leq x \leq 2,~0\leq y\leq 1,~\sqrt{(x-1.4)^2+(y-0.5)^2}\geq 0.18\},\\
			\Omega_3=&\{(x,y)|1.2\leq x \leq 2,~0\leq y\leq 1,~\sqrt{(x-1.4)^2+(y-0.5)^2}< 0.18\}.\\
		\end{split}
	\end{align*}
	The initial solution is given by
	\begin{align*}
		\begin{split}
			&(\rho,u_x,u_y,u_z,B_x,B_y,B_z,p) \left(x,y\right)\\
			&=
			\begin{cases}
				(3.88968,0,0,-0.05234,1,0,3.9353,14.2641),\ & \boldsymbol{x} \in \Omega_1,\\
				(1,-3.3156,0,0,1,0,1,0.04),\ & \boldsymbol{x} \in \Omega_2,\\
				(5,-3.3156,0,0,1,0,1,0.04),\ &\boldsymbol{x} \in \Omega_3.\\
			\end{cases}
		\end{split}
	\end{align*}
	The computational domain is $[0,2]\times[0,1]$, with outflow boundary conditions.
\end{example}

\begin{figure}[htbp!]
	\centering
	\includegraphics[scale=0.35]{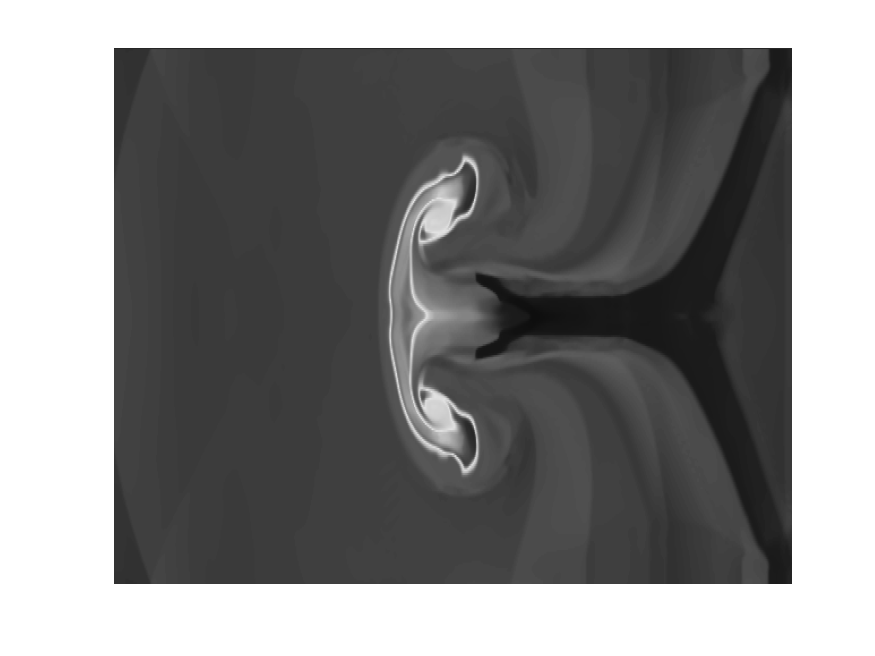}
	\includegraphics[scale=0.35]{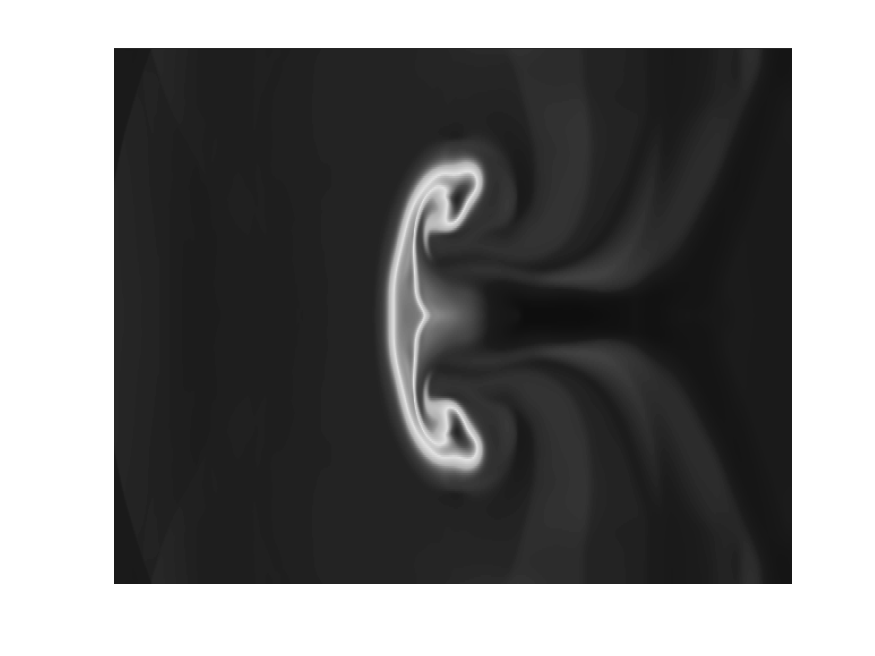}\\
	\includegraphics[scale=0.35]{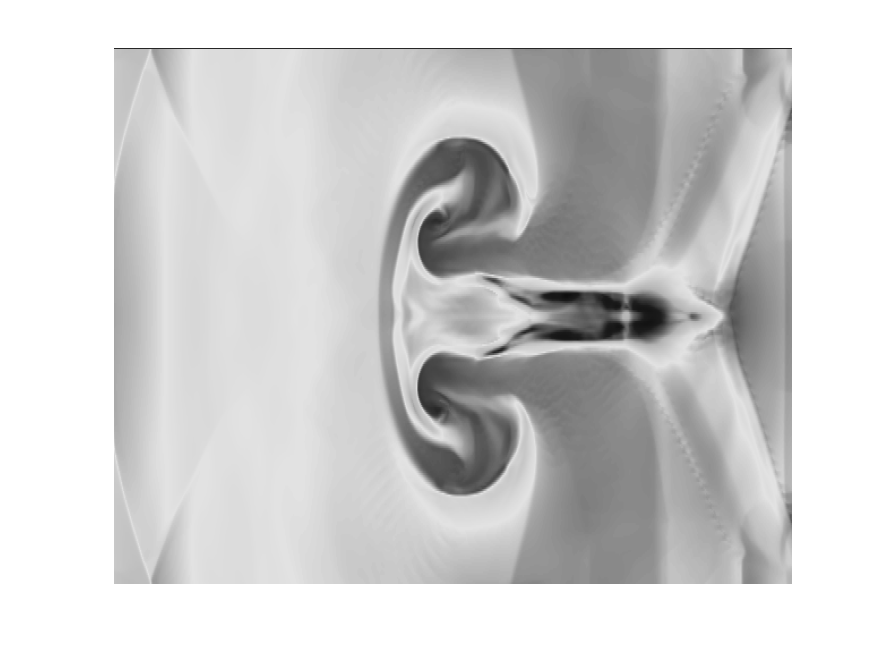}
	\includegraphics[scale=0.35]{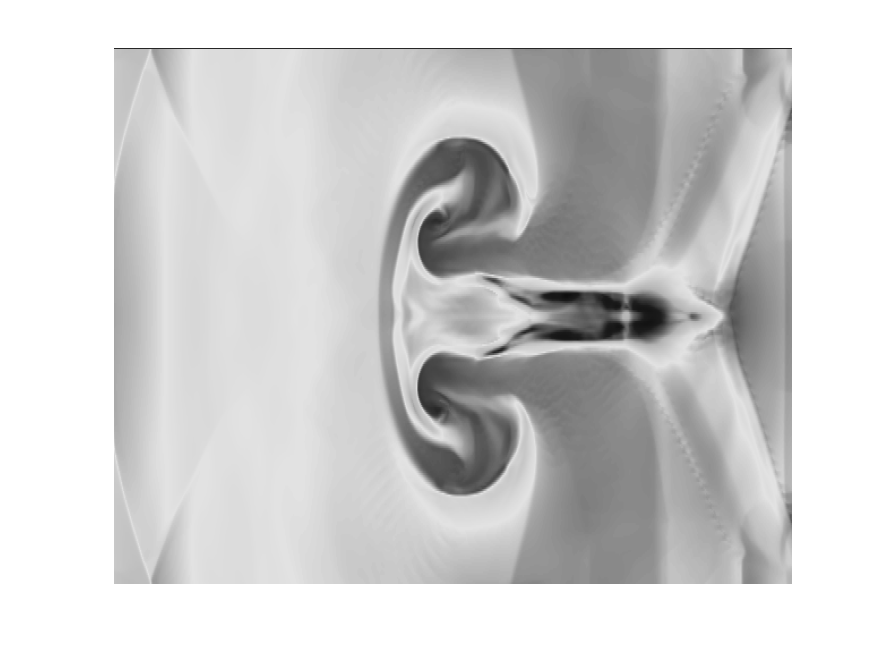}\\
	\includegraphics[scale=0.35]{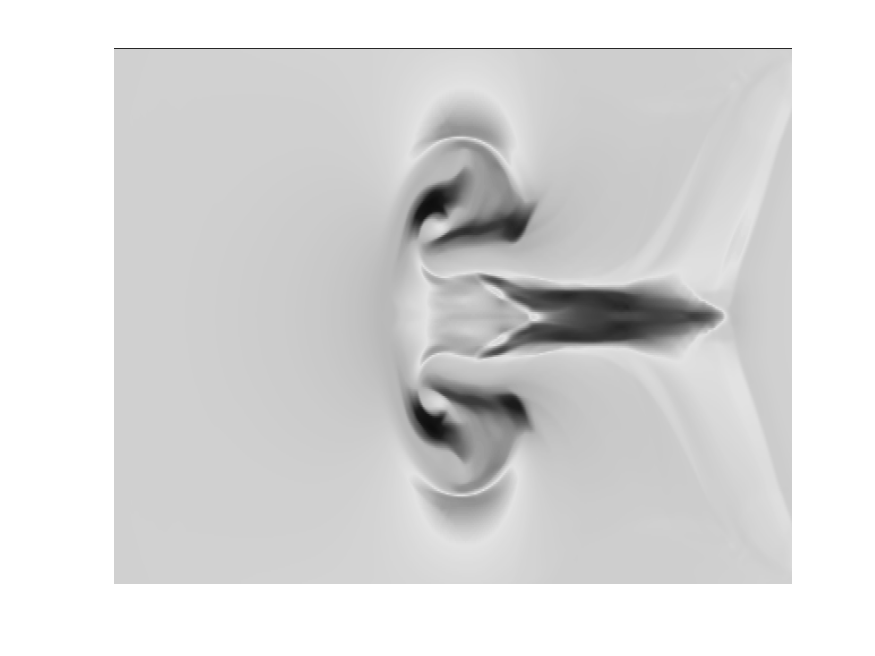}
	\includegraphics[scale=0.35]{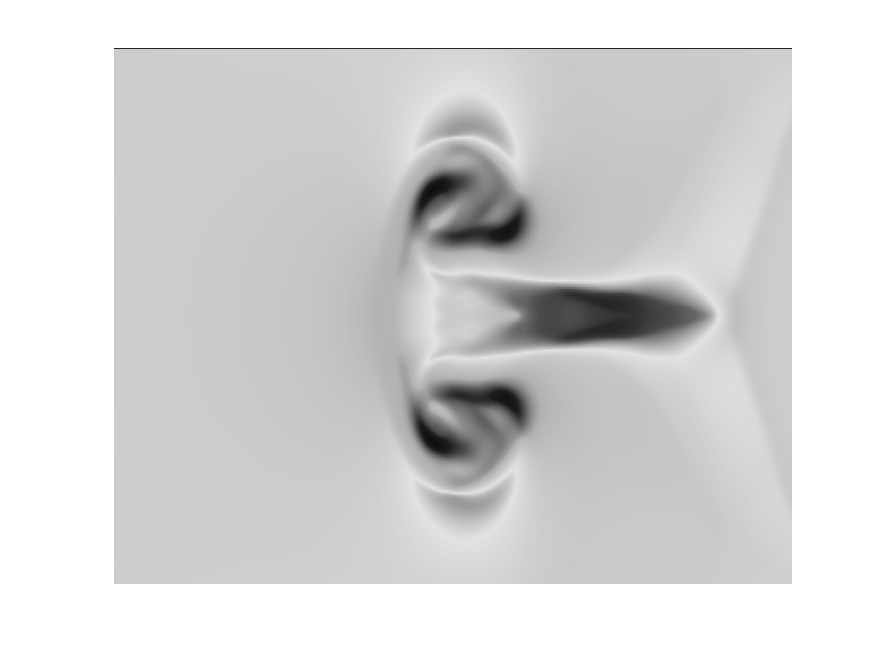}\\
	\includegraphics[scale=0.35]{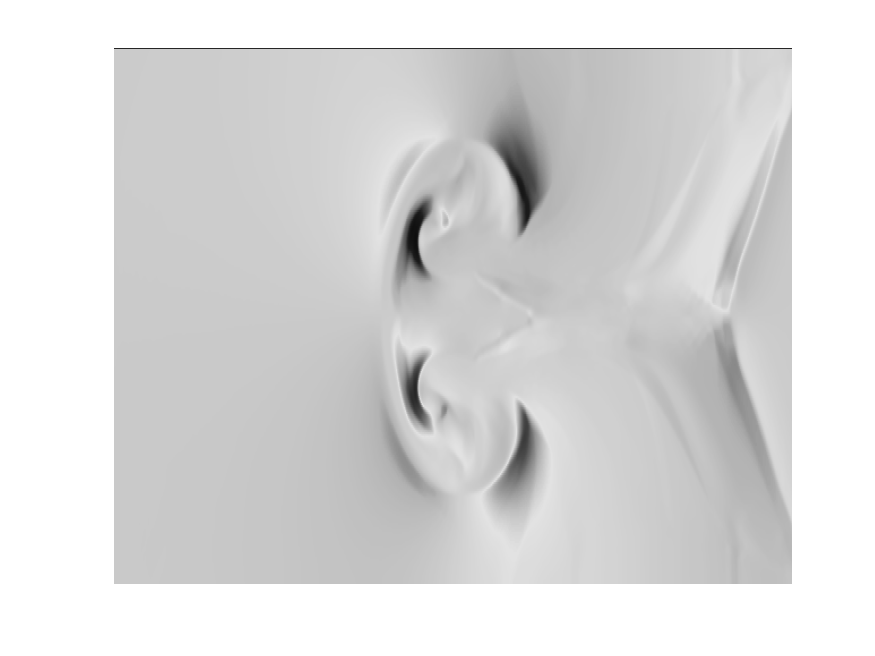}
	\includegraphics[scale=0.35]{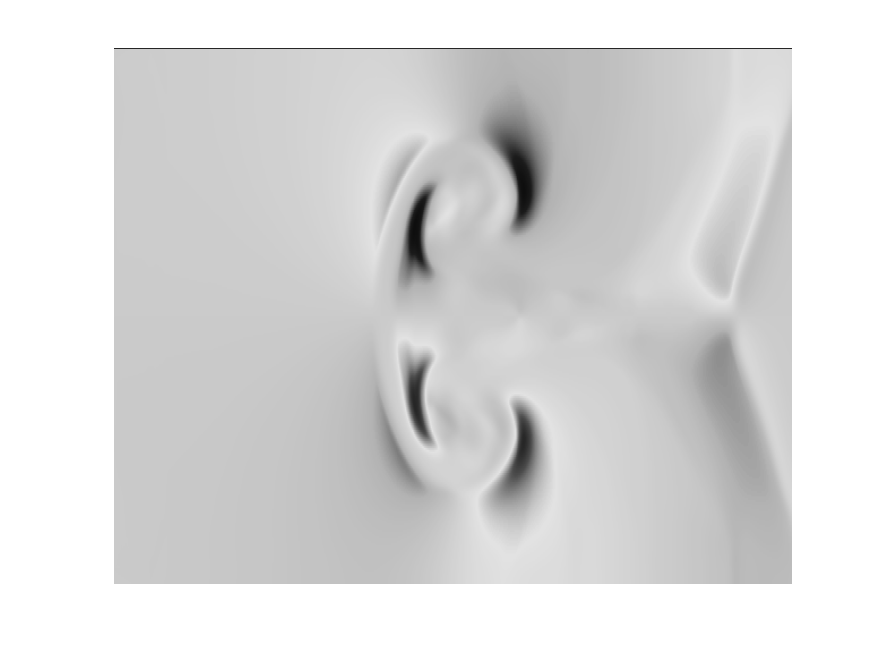}\\
	\caption{Shock cloud interaction problem. From top to bottom: gray-scale images of density $\rho$, pressure $p$, magnetic field components $B_x$ and $B_y$  at $t = 0.6$ on a $600 \times 300$ mesh. Left: LDF-DG with TVB limiter; right: LDF-OEDG.}\label{CloudPlot}
\end{figure}

Simulations are performed on a $600 \times 300$ mesh, up to $t=0.6$. The gray-scale images of density $\rho$, pressure $p$, magnetic field components $B_x$ and $B_y$, computed by the LDF-OEDG method and the LDF-DG method with TVB limiter, are shown in Figure \ref{CloudPlot}. It is observed that the complex flow structures and discontinuities are well resolved, and these results are fairly close to those presented in \cite{Dai1998A,Fu2018Globally}.

\section{Conclusions} \label{Sect5}
This paper presents a locally divergence-free oscillation-eliminating discontinuous Galerkin (LDF-OEDG) method for ideal compressible MHD equations. The OE procedure, which introduces a damping mechanism for modal coefficients of the DG solution, is employed to suppress spurious oscillations near discontinuities. The damping ODE can be exactly solved, thus making the OEDG method stable under normal CFL conditions. The magnetic divergence-free constraint is satisfied by performing a projection of the OEDG solution onto a divergence-free subspace on each element. In each Runge-Kutta stage, the DG spatial discretization, the OE procedure and the LDF projection are fully decoupled, thus enabling an easy implementation of the LDF-OEDG method.  The LDF-OEDG method is applied to a series of benchmark ideal MHD test cases. Numerical results demonstrate the high-order accuracy, strong shock capturing capability and robustness of the LDF-OEDG method. 

\bigskip

\noindent 
{\footnotesize \textbf{Funding} \ Research of Q. Wang is supported by NSFC Grant Nos. 12372284 and U2230402. Research of W. Zeng is supported by the Postdoctoral Fellowship Program of CPSF Grant No. GZC20230213.}

\medskip

\noindent 
{\footnotesize \textbf{Data availability} \ Data will be made available on reasonable request.}

\section*{Declarations}
{\footnotesize \textbf{Conflict of interest} \ The authors have not disclosed any competing interests.}


%
%

\bibliographystyle{spmpsci}      
\bibliography{refs}   

%
%

\end{document}